\documentclass[12pt,a4paper,twoside]{article}

\usepackage{a4wide}
\usepackage[utf8]{inputenc}
\usepackage[T1]{fontenc}
\usepackage{amsmath,amsthm,amsfonts,bbm}
\usepackage{graphicx}
\usepackage{algorithm}
\usepackage{algorithmic}

\newcommand{\R}{\mathbb{R}}
\newcommand{\C}{\mathbb{C}}
\newcommand{\N}{\mathbb{N}}
\newcommand{\abs}[1]{\left\lvert{#1}\right\rvert}
\newcommand{\norm}[1]{{\left\lVert{#1}\right\rVert}}
\newcommand{\spf}[2]{{\left\langle{#1},{#2}\right\rangle}}
\renewcommand{\vec}[1]{\underline{#1}}
\newtheorem{thm}{Theorem}[section]

\begin{document}

\title{Efficient Direct Space-Time Finite Element Solvers for Parabolic Initial-Boundary Value Problems in Anisotropic Sobolev Spaces}
\author{Ulrich~Langer$^1$, Marco~Zank$^2$}
\date{
        $^1$Johann Radon Institute for Computational and Applied Mathematics, \\
        Austrian Academy of Sciences, \\
        Altenberger Straße 69, 4040 Linz, Austria \\[1mm]
        {\tt ulrich.langer@ricam.oeaw.ac.at}  \\[3mm]
        $^2$Fakult\"at f\"ur Mathematik, Universit\"at Wien, \\
        Oskar-Morgenstern-Platz 1, 1090 Wien, Austria \\[1mm]
        {\tt marco.zank@univie.ac.at}
      }
      
\maketitle

\begin{abstract}
    We consider a space-time variational formulation of parabolic initial-boundary value problems in anisotropic Sobolev spaces in combination with a Hilbert-type transformation. This variational setting is the starting point for the space-time Galerkin finite element discretization that leads to a large global linear system of algebraic equations. We propose and investigate new efficient direct solvers for this system. In particular, we use a tensor-product approach with piecewise polynomial, globally continuous ansatz and test functions. The developed solvers are based on the Bartels-Stewart method and on the Fast Diagonalization method, which result in solving a sequence of spatial subproblems. The solver based on the Fast Diagonalization method allows to solve these spatial subproblems in parallel leading to a full parallelization in time. We analyze the complexity of the proposed algorithms, and give numerical examples for a two-dimensional spatial domain, where sparse direct solvers for the spatial subproblems are used.
\end{abstract}

\section{Introduction}

Parabolic initial-boundary value problems are usually discretized by time-stepping schemes and spatial finite element methods. These methods treat the time and spatial variables differently, see, e.g., \cite{Thomee2006}. In addition, the resulting approximation methods are sequential in time. In contrast to these approaches, space-time methods discretize time-dependent partial differential equations without separating the temporal and spatial directions. In particular, they are based on space-time variational formulations. There exist various space-time techniques for parabolic problems, which are based on variational formulations in Bochner-Sobolev spaces, see, e.g., \cite{Andreev2013, LZ:BankVassilevskiZikatanov:2016a, fuehrer2019spacetime, LZ:HughesFrancaHulbert:1989a, Larsson2017, LoliMontardiniSangalliTani2019, Mollet2014, SchwabStevenson2009, Steinbach2015, Stevenson2020, UrbanPatera2014}, or on discontinuous Galerkin methods, see, e.g., \cite{HoferLangerNeumuellerSchneckenleitner2019, NeumuellerDissBuch2013, NeumuellerSmears2019}, or on discontinuous Petrov-Galerkin methods, see, e.g., \cite{LZ:FuehrerHeuerSenGupta:2017a}, and the references therein. We refer the reader to \cite{LZ:Gander:2015a} and \cite{LZ:SteinbachYang:2019a} for a comprehensive overview of parallel-in-time and space-time methods, respectively. An alternative is the discretization of space-time variational formulations in anisotropic Sobolev spaces, see, e.g., \cite{Devaud2019, LarssonSchwab2015, SchwabStevenson2017, SteinbachZank2020Coercive, ZankDissBuch2020}. These variational formulations allow the complete analysis of inhomogeneous Dirichlet or Neumann conditions, and were used for the analysis of the resulting boundary integral operators, see \cite{Costabel1990, DohrSteinbachKazuki2019}. Hence, discretizations for variational formulations in anisotropic Sobolev spaces can be used for the interior problems of FEM-BEM couplings for transmission problems.

In this work, the approach in anisotropic Sobolev spaces is applied in combination with a novel Hilbert-type transformation operator $\mathcal H_T$, which has recently been introduced in \cite{SteinbachZank2020Coercive, ZankDissBuch2020}. This transformation operator $\mathcal H_T$ maps the ansatz space to the test space, and gives a symmetric and elliptic variational setting of the first-order time derivative. The homogeneous Dirichlet problem for the nonstationary diffusion respectively heat equation
\begin{equation} \label{Einf:WaermePDG}
    \left. \begin{array}{rclcl}
    \partial_t u(x,t) - \Delta_xu(x,t) & = & f(x,t) & \quad & \mbox{for} \;
    (x,t) \in Q = \Omega \times (0,T), \\[1mm]
    u(x,t) & = & 0 & & \mbox{for} \; (x,t) \in \Sigma = \partial\Omega \times [0,T], 
    \\[1mm] 
    u(x,0) & = & 0 & & \mbox{for} \; x \in \Omega,
    \end{array} \right \}
\end{equation}
serves as model problem for a parabolic initial-boundary value problem, where $\Omega \subset \R^d,$ $d=1,2,3,$ is a bounded Lipschitz domain with boundary $\partial\Omega$, $T>0$ is a given terminal time, and $f$ is a given right-hand side. With the help of the Hilbert-type transformation operator $\mathcal H_T$, a Galerkin finite element method is derived, which results in one global linear system
\begin{equation} \label{Einf:LGS}
    K_h \vec u = \vec f.
\end{equation}
When using a tensor-product approach, the system matrix $K_h$ can be represented as a sum of Kronecker products. The purpose of this paper is the development of efficient direct space-time solvers for the global linear system \ref{Einf:LGS}, exploiting the Kronecker structure of $K_h$. Therefore, we apply the Bartels-Stewart method \cite{BartelsStewart1972} and the Fast Diagonalization method \cite{Lynch1964} to solve \eqref{Einf:LGS}, see also \cite{Gardiner1992, Simoncini2016}. 
For both methods, we derive complexity estimates of the resulting algorithms.

The rest of this paper is organized as follows: In Section~\ref{Sec:VF}, we consider the space-time variational formulation in anisotropic Sobolev spaces and the Hilbert-type transformation operator $\mathcal H_T$ with its main properties. In Section~\ref{Sec:Vor}, we rephrase properties of the Kronecker product and of sparse direct solvers, which are needed for the new space-time solver. Section~\ref{Sec:Loeser} is devoted to the construction of efficient space-time solvers. Numerical examples for a two-dimensional spatial domain are presented in Section~\ref{Sec:Num}. Finally, we draw some conclusions in Section~\ref{Sec:Zum}.

\section{Space-Time Method in Anisotropic Sobolev Spaces} \label{Sec:VF}

In this section, we give the variational setting for the parabolic model problem \eqref{Einf:WaermePDG}, which is studied in greater detail in \cite{Costabel1990, Ladyzhenskaya1968, LM11968, LM21968, SteinbachZank2020Coercive, ZankDissBuch2020}. We consider the space-time variational formulation of \eqref{Einf:WaermePDG} in anisotropic Sobolev spaces to find $u \in H^{1,1/2}_{0;0,\,}(Q)$ such that 
\begin{equation}\label{VF:VF}
	a(u,v) = \spf{f}{v}_{Q}
\end{equation}
for all $v \in H^{1,1/2}_{0;\,,0}(Q),$ where $f \in [H^{1,1/2}_{0;\,,0}(Q)]'$ is a given right-hand side. Here, the bilinear form $a(\cdot,\cdot) \colon \, H^{1,1/2}_{0;0,\,}(Q) \times H^{1,1/2}_{0;\,,0}(Q) \to \R,$
\begin{equation*}
	a(u,v) := \spf{\partial_{t}u}{v}_{Q} + \spf{\nabla_x u}{\nabla_x v}_{L^2(Q)} 
\end{equation*}
for $u \in H^{1,1/2}_{0;0,\,}(Q)$, $v \in H^{1,1/2}_{0;\,,0}(Q)$, is bounded, i.e. there exists a constant $C>0$ such that
\begin{equation*}
	\forall u \in H^{1,1/2}_{0;0,\,}(Q) \colon \, \forall v \in H^{1,1/2}_{0;\,,0}(Q)\colon \quad \abs{a(u,v)} \leq C \norm{u}_{H^{1,1/2}_{0;0,\,}(Q)} \norm{v}_{H^{1,1/2}_{0;\,,0}(Q)},
\end{equation*}
see \cite[Lemma 2.6, p. 505]{Costabel1990}. The anisotropic Sobolev spaces
\begin{align*}
	H^{1,1/2}_{0;0,\,}(Q) &:= H^{1/2}_{0,} (0,T; L^2(\Omega)) \cap L^2(0,T; H^1_0(\Omega)), \\
	H^{1,1/2}_{0;\,,0}(Q) &:= H^{1/2}_{,0} (0,T; L^2(\Omega)) \cap L^2(0,T; H^1_0(\Omega)) 
\end{align*}
are endowed with the Hilbertian norms
\begin{align*}
	\| v \|_{ H^{1,1/2}_{0;0,\,}(Q) } &:= \sqrt{ \| v \|_{H^{1/2}_{0,}(0,T; L^2(\Omega))}^2 + \norm{\nabla_x v}_{L^2(Q)}^2 }, \\
	\| w \|_{ H^{1,1/2}_{0;\,,0}(Q) } &:= \sqrt{ \| w \|_{H^{1/2}_{,0}(0,T; L^2(\Omega))}^2 + \norm{\nabla_x w}_{L^2(Q)}^2 }
\end{align*}
with the usual Bochner-Sobolev norms
\begin{align*}
	\| v \|_{H^{1/2}_{0,}(0,T; L^2(\Omega))} &:= \sqrt{\norm{v}_{H^{1/2}(0,T; L^2(\Omega))}^2 + \int_0^T \frac{\norm{v(\cdot, t)}_{L^2(\Omega)}^2}{t} \mathrm dt}, \\
        \| w \|_{H^{1/2}_{,0}(0,T; L^2(\Omega))} &:= \sqrt{\norm{w}_{H^{1/2}(0,T; L^2(\Omega))}^2 + \int_0^T \frac{\norm{w(\cdot,t)}_{ L^2(\Omega)}^2}{T-t} \mathrm dt},
\end{align*}
see \cite{LM11968, LM21968, SteinbachZank2020Coercive, ZankDissBuch2020} for more details. The dual space $[H^{1,1/2}_{0;\,,0}(Q)]'$ is characterized as completion of $L^2(Q)$ with respect to the Hilbertian norm
\begin{equation*}
	\| f \|_{[H^{1,1/2}_{0;\,,0}(Q)]'} := \sup_{0 \neq w \in H^{1,1/2}_{0;\,,0}(Q) } 
\frac{\abs{\langle f, w \rangle_Q}}{ \| w \|_{ H^{1,1/2}_{0;\,,0}(Q) } },
\end{equation*}
where $\langle \cdot , \cdot \rangle_Q$ denotes the duality pairing as extension of the inner product in $L^2(Q).$ In \cite{Costabel1990}, the following existence and uniqueness theorem is proven by a transposition and interpolation argument as in \cite{LM11968, LM21968}, see also \cite{Ladyzhenskaya1968}.

\begin{thm}
	Let the right-hand side $f \in [H^{1,1/2}_{0;\,,0}(Q)]'$ be given. Then, the variational formulation \eqref{VF:VF} has a unique solution $u \in H^{1,1/2}_{0;0,\,}(Q),$ satisfying
  \begin{equation*}
    \norm{u}_{H^{1,1/2}_{0;0,\,}(Q)} \leq C \norm{f}_{[H^{1,1/2}_{0;\,,0}(Q)]'} 
  \end{equation*}
  with a constant $C>0.$ Furthermore, the solution operator
  \begin{equation*}
		 \mathcal L \colon \, [H^{1,1/2}_{0;\,,0}(Q)]' \to H^{1,1/2}_{0;0,\,}(Q), \quad \mathcal L f:=u,
  \end{equation*}
  is an isomorphism.
\end{thm}

For simplicity, we only consider homogeneous Dirichlet conditions, where inhomogeneous Dirichlet conditions can be treated via homogenization as for the elliptic case, see \cite[p. 61-62]{Steinbach2008}, since for any Dirichlet data $g \in H^{1/2,1/4}(\Sigma)$, an extension $u_g \in H^{1,1/2}_{\,; 0,\,}(Q)$ with $u_{g|\Sigma} = g$ exists, see \cite{Costabel1990, DohrSteinbachKazuki2019} for more details.

For a discretization scheme, let the bounded Lipschitz domain $\Omega \subset \R^d$ be an interval $\Omega = (0,L)$ for $d=1,$ or polygonal for $d=2,$ or polyhedral for $d=3.$ For a tensor-product ansatz, we consider admissible decompositions
\begin{equation*}
    \overline{Q} = \overline{\Omega} \times [0,T] = \bigcup_{i=1}^{N_x}\overline{\omega_i} \times \bigcup_{\ell=1}^{N_t} [t_{\ell-1},t_\ell] 
\end{equation*}
with $N:=N_x \cdot N_t$ space-time elements, where the time intervals $(t_{\ell-1},t_\ell)$ with mesh sizes $h_{t,\ell} = t_\ell - t_{\ell-1}$ are defined via the decomposition
\begin{equation*}
    0 = t_0 < t_1 < t_2 < \dots < t_{N_t -1} < t_{N_t} = T
\end{equation*}
of the time interval $(0,T)$. The maximal and the minimal time mesh sizes are denoted by $h_t := h_{t,\max} := \max_{\ell} h_{t,\ell}$ and $h_{t,\min} := \min_{\ell} h_{t,\ell}$, respectively. For the spatial domain $\Omega$, 
we consider a shape-regular sequence $(\mathcal T_\nu)_{\nu \in {\mathbb{N}}}$ of admissible decompositions
\begin{equation*}
    \mathcal T_\nu := \{ \omega_i \subset \R^{d} \colon i=1,\dots,N_x\}
\end{equation*}
of $\Omega$ into finite elements $\omega_i \subset \R^d$ with mesh sizes $h_{x,i}$
and the maximal mesh size $h_x := \max_{i} h_{x,i}$.
The spatial elements $\omega_{i}$ are 
intervals for $d=1$, triangles or quadrilaterals for $d=2$, and tetrahedra 
or hexahedra for $d=3$. 
Next, we introduce the finite element space
\begin{equation} \label{VF:Q1}
 Q_h^1(Q) := V_{h_x,0}^1(\Omega) \otimes S_{h_t}^1(0,T)
\end{equation}
of piecewise multilinear, continuous functions, i.e.
\begin{equation*}
    V_{h_x,0}^1(\Omega) = \mbox{span} \{ \psi_j^1 \}_{j=1}^{M_x} \subset H^1_0(\Omega), \quad S_{h_t}^1(0,T) = \mbox{span} \{ \varphi_\ell^1 \}_{\ell=0}^{N_t} \subset H^1(0,T).
\end{equation*}
In fact, $V_{h_x,0}^1(\Omega)$ is either the space $S_{h_x}^1(\Omega) \cap H^1_0(\Omega)$ of piecewise 
linear, continuous functions on intervals ($d=1$), triangles 
($d=2$), and tetrahedra ($d=3$), or $V_{h_x,0}^1(\Omega)$ is the space 
$Q_{h_x}^1(\Omega) \cap H^1_0(\Omega)$ of piecewise linear/bilinear/trilinear, continuous functions on intervals ($d=1$), 
quadrilaterals ($d=2$), and hexahedra ($d=3$). Analogously, for a fixed polynomial degree $p \in \N$, we consider the space of piecewise polynomial, continuous functions
\begin{equation} \label{VF:Qh_p}
    Q_h^p(Q) := V_{h_x,0}^p(\Omega) \otimes S_{h_t}^p(0,T).
\end{equation}
Using the finite element space \eqref{VF:Q1}, it turns out that a discretization of \eqref{VF:VF} with the conforming ansatz space $Q_h^1(Q) \cap H^{1,1/2}_{0;0,\,}(Q)$ and the conforming test space $Q_h^1(Q) \cap H^{1,1/2}_{0;\,,0}(Q)$ is not stable, see \cite[Section 3.3]{ZankDissBuch2020}. A possible way out is the modified Hilbert transformation ${\mathcal{H}}_T$ defined by
\begin{equation*}
	({\mathcal{H}}_Tu)(x,t) :=  \sum_{i=1}^\infty \sum_{k=0}^\infty u_{i,k} \cos \left( \Big( \frac{\pi}{2} + k\pi \Big) \frac{t}{T} \right) \phi_i(x), \quad (x,t) \in Q,
\end{equation*}
where the given function $u \in L^2(Q)$ is represented by
\begin{equation} \label{VF:L2_Darstellung}
	u(x,t) = \sum_{i=1}^\infty \sum_{k=0}^\infty u_{i,k} \sin \left( \Big( \frac{\pi}{2} + k\pi \Big) \frac{t}{T} \right) \phi_i(x), \quad (x,t) \in Q,
\end{equation}
with the eigenfunctions $\phi_i \in H^1_0(\Omega)$ and eigenvalues $\mu_i \in \R$, satisfying 
\begin{equation*}
- \Delta \phi_i = \mu_i \phi_i  \quad \text{ in } \Omega, \quad \phi_i = 0 \quad \text{ on } \partial \Omega, \quad \norm{\phi_i}_{L^2(\Omega)} = 1, \quad i \in \N.
\end{equation*}
This approach was introduced recently in \cite{SteinbachZank2020Coercive} and \cite[Section 3.4]{ZankDissBuch2020}. The novel transformation $\mathcal H_T$ acts on the finite terminal $(0,T)$, whereas analogous considerations of an infinite time interval $(0,\infty)$ with the classical Hilbert transformation are investigated in \cite{Devaud2019, Fontes1996, Fontes2009, LarssonSchwab2015}. The most important properties of $\mathcal H_T$ are summarized in the following, see \cite{SteinbachZank2020Coercive, SteinbachZank2020Note, Zank2020Exact, ZankDissBuch2020}. The map
\begin{equation*}
	{\mathcal{H}}_T \colon \, H^{1,1/2}_{0;0,}(Q) \to H^{1,1/2}_{0;\,,0}(Q)
\end{equation*}
is norm preserving, bijective and fulfills the coercivity property
\begin{equation} \label{VF:HT_Abl_posdef}
    \spf{\partial_t u}{\mathcal{H}_T v}_{Q} = \frac{1}{2} \sum_{i=1}^\infty \sum\limits_{k=0}^\infty \Big( \frac{\pi}{2} + k\pi \Big) u_{i,k} \cdot v_{i,k}=:\spf{u}{v}_{H^{1/2}_{0,}(0,T; L^2(\Omega)),F}
\end{equation}
for functions $u,v \in H^{1/2}_{0,}(0,T; L^2(\Omega))$ with expansion coefficients $u_{i,k}, v_{i,k}$ as in \eqref{VF:L2_Darstellung}. Note that the norm induced by the inner product $\spf{\cdot}{\cdot}_{H^{1/2}_{0,}(0,T; L^2(\Omega)),F}$ is equivalent to the norm $\norm{\cdot}_{H^{1/2}_{0,}(0,T; L^2(\Omega))}.$
Moreover, the relations
\begin{align} 
  \forall v \in L^2(Q) \colon& \quad \spf{v}{\mathcal{H}_T v}_{L^2(Q)} \geq 0, \nonumber \\
  \forall s > 0 \colon \forall v \in H^s_{0,}(0,T; L^2(\Omega)), v \neq 0 \colon& \quad \spf{v}{\mathcal{H}_T v}_{L^2(Q)} > 0  \label{VF:HT_posdef}
\end{align}
hold true. With the modified Hilbert transformation ${\mathcal{H}}_T$, the variational formulation \eqref{VF:VF} is equivalent to find $u \in H^{1,1/2}_{0;0,\,}(Q)$ such that 
\begin{equation}\label{VF:VF_HT}
	\forall v \in H^{1,1/2}_{0;0,\,}(Q) \colon \, a(u, \mathcal{H}_T v) = \spf{f}{{\mathcal{H}}_T v}_Q. 
\end{equation}
Hence, unique solvability of the variational formulation \eqref{VF:VF_HT} follows from the unique solvability of \eqref{VF:VF}, which implies the stability estimate
\begin{equation*}
	\forall  u \in H^{1,1/2}_{0;0,}(Q)  \colon \, c \, \| u \|_{H^{1,1/2}_{0;0,}(Q)} \leq \sup\limits_{0 \neq v \in H^{1,1/2}_{0;0,}(Q)} \frac{ \abs{a(u, \mathcal{H}_T v) }}{\| v \|_{H^{1,1/2}_{0;0,}(Q)}}
\end{equation*}
with a constant $c>0.$ When using some conforming space-time finite element space ${\mathcal{V}}_h \subset  H^{1,1/2}_{0;0,}(Q),$ the Galerkin variational formulation of \eqref{VF:VF_HT} is to find $u_h \in {\mathcal{V}}_h$ such that
\begin{equation}\label{VF:VF_HT_disk}
    \forall v_h \in {\mathcal{V}}_h \colon \, a(u_h,\mathcal{H}_T v_h) = \spf{f}{\mathcal{H}_T v_h}_Q.
\end{equation}
Note that ansatz and test spaces are equal. In \cite{ZankDissBuch2020}, the following theorem is proven.
\begin{thm}\label{VF:Thm:BelDisk}
	Let ${\mathcal{V}}_h \subset  H^{1,1/2}_{0;0,}(Q)$ be a conforming space-time finite element space and let $f \in [H^{1,1/2}_{0; \, ,0}(Q)]'$ be a given right-hand side. Then, a unique solution $u_h \in {\mathcal{V}}_h$ of the Galerkin variational formulation \eqref{VF:VF_HT_disk} exists. If, in addition, the right-hand side fulfills $f \in [H^{1/2}_{,0}(0,T;L^2(\Omega))]' \subset [H^{1,1/2}_{0; \, ,0}(Q)]',$ then the stability estimate
	\begin{equation*}
		\| u_h \|_{H^{1/2}_{0,}(0,T;L^2(\Omega))} \leq c \| f \|_{[H^{1/2}_{,0}(0,T;L^2(\Omega))]'}
	\end{equation*}
	is true with a constant $c>0.$
\end{thm}
Theorem~\ref{VF:Thm:BelDisk} states that, under the assumption $f \in [H^{1/2}_{,0}(0,T;L^2(\Omega))]'$, any conforming space-time finite element space ${\mathcal{V}}_h \subset  H^{1,1/2}_{0;0,}(Q)$ leads to an unconditionally stable method, i.e. no CFL condition is required. For the choice of the tensor-product space-time finite element space
\begin{equation*}
 {\mathcal{V}}_h = Q_h^p(Q) \cap H^{1,1/2}_{0;0,\,}(Q)
\end{equation*}
from \eqref{VF:Qh_p}, the Galerkin variational formulation \eqref{VF:VF_HT_disk} to find $u_h \in  Q_h^p(Q) \cap H^{1,1/2}_{0;0,\,}(Q)$ such that
\begin{equation*}
    \forall v_h \in  Q_h^p(Q) \cap H^{1,1/2}_{0;0,\,}(Q) \colon \, a(u_h,\mathcal{H}_T v_h) = \spf{f}{\mathcal{H}_T v_h}_Q
\end{equation*}
fulfills the space-time error estimates
  \begin{align}
    \norm{u - u_h}_{H^{1/2}_{0,}(0,T;L^2(\Omega))} &\leq c \, h^{p+1/2}, \label{VF:Fehler_H12} \\
    \norm{u - u_h}_{L^2(Q)} &\leq c \, h^{p+1}, \label{VF:Fehler_L2} \\
    \abs{u - u_h}_{H^1(Q)} &\leq c \, h^{p} \label{VF:Fehler_H1}
  \end{align}
with $h= \max\{h_t, h_x\}$ and with a constant $c>0$ for a sufficiently smooth solution $u \in H^{1,1/2}_{0;0,}(Q)$ of \eqref{VF:VF} and a sufficiently regular boundary $\partial \Omega,$ where for the $H^1(Q)$ error estimate \eqref{VF:Fehler_H1}, the sequence $(\mathcal T_\nu)_{\nu \in {\mathbb{N}}}$ of decompositions of $\Omega$ is additionally assumed to be globally quasi-uniform, see \cite{SteinbachZank2020Coercive, ZankDissBuch2020} for details.

In the remainder of this work, we consider $p=1$, i.e. the tensor-product space of piecewise linear, continuous functions ${\mathcal{V}}_h = Q_h^1(Q) \cap H^{1,1/2}_{0;0,\,}(Q)$, where analogous results hold true for an arbitrary polynomial degree $p>1.$ So, the number of the degrees of freedom is given by
\begin{equation*}
    \mathrm{dof} = N_t \cdot M_x.
\end{equation*}
For an easier implementation, we approximate the right-hand side $f\in L^2(Q)$ by
\begin{equation} \label{VF:Projektion}
    f \approx Q_{h}^0 f = \sum_{j=1}^{N_x} \sum_{\ell=1}^{N_t} f_{j,\ell} \,\psi_j^0 \varphi_\ell^0 \in S_{h_x}^0(\Omega) \otimes S_{h_t}^0(0,T)
\end{equation}
with coefficients $f_{j,\ell} \in \R$, where $Q_{h}^0 \colon \, L^2(Q) \to S_{h_x}^0(\Omega) \otimes S_{h_t}^0(0,T)$ is the $L^2(Q)$ projection on the piecewise constant functions $S_{h_x}^0(\Omega) \otimes S_{h_t}^0(0,T)$ with $S_{h_x}^0(\Omega) = \mbox{span} \{ \psi_j^0 \}_{j=1}^{N_x}$ and $S_{h_t}^0(0,T) = \mbox{span} \{ \varphi_\ell^0 \}_{\ell=1}^{N_t}$. So, we consider the perturbed variational formulation to find $\widetilde u_h \in  Q_h^1(Q) \cap H^{1,1/2}_{0;0,\,}(Q)$ such that
\begin{equation}\label{VF:VF_HT_disk_Tensor_gestoert}
    \forall v_h \in  Q_h^1(Q) \cap H^{1,1/2}_{0;0,\,}(Q) \colon \, a(\widetilde u_h,\mathcal{H}_T v_h) = \spf{Q_{h}^0 f}{\mathcal{H}_T v_h}_{L^2(Q)}.
\end{equation}
Note that, for piecewise linear functions, i.e. $p=1$, the space-time error estimates \eqref{VF:Fehler_H12}, \eqref{VF:Fehler_L2}, \eqref{VF:Fehler_H1} are not spoilt. Note additionally that, for $p>1$, a projection on polynomials of degree $p-1$ should be used instead of $Q_{h}^0$ for preserving the space-time error estimates \eqref{VF:Fehler_H12}, \eqref{VF:Fehler_L2}, \eqref{VF:Fehler_H1}. After an appropriate ordering of the degrees of freedom, the discrete variational formulation \eqref{VF:VF_HT_disk_Tensor_gestoert} is equivalent to the global linear system
\begin{equation} \label{VF:LGS}
  K_h \vec{u} = \vec{\widetilde F}^{\mathcal H_T}
\end{equation}
with the system matrix
\begin{equation*}
  K_h=A_{h_t}^{\mathcal H_T} \otimes M_{h_x} + M_{h_t}^{\mathcal H_T} \otimes A_{h_x}  \in \R^{ N_t \cdot M_x \times N_t \cdot M_x},
\end{equation*}
where $M_{h_x} \in \R^{M_x \times M_x}$ and $A_{h_x} \in \R^{M_x \times M_x}$ denote spatial mass and stiffness matrices given by
\begin{equation} \label{VF:Matritzen_Mass_Stiff}
	M_{h_x}[i,j] = \langle \psi_j^1, \psi_i^1 \rangle_{L^2(\Omega)}, \quad A_{h_x}[i,j] = \langle \nabla_x \psi_j^1, \nabla_x \psi_i^1 \rangle_{L^2(\Omega)}, \quad i,j=1,\dots, M_x,
\end{equation}
and $M_{h_t}^{\mathcal H_T} \in \R^{N_t \times N_t}$ and $A_{h_t}^{\mathcal H_T} \in \R^{N_t \times N_t}$ are defined by
\begin{equation*}
  M_{h_t}^{\mathcal H_T}[\ell, k] := \spf{\varphi_k^1}{\mathcal H_T \varphi_\ell^1}_{L^2(0,T)}, \quad  A_{h_t}^{\mathcal H_T}[\ell, k] := \spf{\partial_t \varphi_k^1}{\mathcal H_T \varphi_\ell^1}_{L^2(0,T)}
\end{equation*}
for $\ell,k=1,\dots,N_t$. Note that the matrix $A_{h_t}^{\mathcal H_T}$ is dense, symmetric and positive definite, see \eqref{VF:HT_Abl_posdef}, whereas the matrix $M_{h_t}^{\mathcal H_T}$ is dense, nonsymmetric and positive definite, see \eqref{VF:HT_posdef}. Additionally, the vector of the right-hand side in \eqref{VF:LGS} is given by
\begin{equation*}
    \vec{\widetilde F}^{\mathcal H_T}  := \left( \vec{\widetilde f}_1, \dots, \vec{\widetilde f}_{M_x}  \right)^\top  \in \R^{N_t \cdot M_x}
\end{equation*}
with the vectors $\vec{\widetilde f}_i \in \R^{N_t}$, \, $i=1, \dots, M_x$, where, with the help of \eqref{VF:Projektion},
\begin{equation*}
    \vec{\widetilde f}_i[k] := \spf{Q_{h}^0 f}{\psi_i^1\, \mathcal{H}_T \varphi_k^1}_{L^2(Q)} = \sum_{j=1}^{N_x} \sum_{\ell=1}^{N_t} f_{j,\ell} \spf{\psi_j^0}{\psi_i^1}_{L^2(\Omega)} \spf{\varphi_\ell^0}{\mathcal{H}_T \varphi_k^1}_{L^2(0,T)},
\end{equation*}
$k=1,\dots, N_t.$ To assemble the vector of the right-hand side in \eqref{VF:LGS}, the relation
\begin{equation*}
    \vec{\widetilde f}_i[k] := \widetilde F[i,k], \qquad i=1, \dots, M_x, \, k=1,\dots, N_t,
\end{equation*}
holds true with $\widetilde{F} := M_{h_x}^{1,0} F (C_{h_t}^{\mathcal H_T})^\top \in \R^{M_x \times N_t},$ where
\begin{align*}
  M_{h_x}^{1,0}[i,j] &:= \spf{\psi_j^0}{\psi_i^1}_{L^2(\Omega)}, && i=1,\dots, M_x, \, j=1,\dots,N_x, \\
  F[j,\ell] &:= f_{j,\ell}, && j=1,\dots,N_x, \,\ell=1,\dots, N_t, \\
  C_{h_t}^{\mathcal H_T}[k,\ell] &:= \spf{\varphi_\ell^0}{\mathcal{H}_T \varphi_k^1}_{L^2(0,T)}, && k=1,\dots,N_t, \,\ell=1,\dots, N_t.
\end{align*}

\section{Preliminaries for the Space-Time Solvers}\label{Sec:Vor}

In this section, some properties of the Kronecker product and direct solvers, which are needed in Section~\ref{Sec:Loeser}, are summarized.

\subsection{Kronecker Product}

In this subsection, some basic properties of the Kronecker product are stated, see, e.g., \cite{HardySteeb2019, Simoncini2016}. Let $A,C \in \C^{N_A \times N_A}$, \,$B,D \in \C^{N_B \times N_B}$ and $X \in \C^{N_A \times N_B}$ be given matrices for $N_A, N_B \in \N.$ The Kronecker product is defined as the matrix
\begin{equation*}
 A \otimes B := \begin{pmatrix}
                    A[1,1] B & A[1,2] B & \cdots & A[1,N_A] B \\
                    A[2,1] B & A[2,2] B & \cdots & A[2,N_A] B \\
                    \vdots   &  \vdots  & \ddots & \vdots     \\
                    A[N_A,1] B & A[N_A,2] B & \cdots & A[N_A,N_A] B
                \end{pmatrix}
                \in \C^{N_A \cdot N_B \times N_A \cdot N_B}.
\end{equation*}
Furthermore, the vectorization of a matrix converts the matrix into a column vector, i.e. we define
\begin{equation*}
 \mathrm{vec}(X) := (X[1,1], X[2,1], \dots, X[N_A,1], X[1,2], \dots, X[N_A,N_B])^\top \in \C^{N_A \cdot N_B \times 1}.
\end{equation*}
In the remainder of this work, we use the following properties of the Kronecker product and the vectorization of a matrix:
\begin{itemize}
    \item For the conjugate transposition and transposition, it holds true that
    \begin{equation*}
        (A \otimes B)^* = A^* \otimes B^* \quad \text{ and } \quad (A \otimes B)^\top = A^\top \otimes B^\top.
    \end{equation*}
    \item For regular matrices $A, B$, we have
    \begin{equation*}
      (A \otimes B)^{-1} = A^{-1} \otimes B^{-1}.
    \end{equation*}
    \item The mixed-product property
    \begin{equation*}
       (A \otimes B) (C \otimes D) = (AC) \otimes (BD)
    \end{equation*}
    is valid.
    \item It holds true that
    \begin{equation} \label{Vor:Vektorisierung}
        \mathrm{vec}(A X B) = (B^\top \otimes A) \mathrm{vec}(X),
    \end{equation}
    where also in the case of complex matrices, only the transposition is applied.
\end{itemize}
For a given vector $\vec v \in \C^{N_A \cdot N_B}$, define the matrix
\begin{equation*}
    V:= \begin{pmatrix}
          \vec v_1 & \vec v_2 & \cdots & \vec v_{N_B}
        \end{pmatrix}
        \in \C^{N_A \times N_B}
\end{equation*}
with $\vec v_i \in \C^{N_A}$ given by $\vec v_i[k] = \vec v[(i-1)N_A + k]$ for $k=1,\dots,N_A,$ $i=1,\dots,N_B$, i.e. $\mathrm{vec}(V) = \vec v.$ Then, the equality \eqref{Vor:Vektorisierung} yields 
\begin{equation} \label{Vor:VektorisierungVektor}
 (B^\top \otimes A) \vec v = (B^\top \otimes A) \mathrm{vec}(V) = \mathrm{vec}(A V B).
\end{equation}

\subsection{Sparse Direct Solver} \label{Sec:Vor:SparseLoeser}

In this subsection, we repeat some properties of sparse direct solver, like left-looking/right-looking/multifrontal methods, for solving linear systems $K_{h_x} \vec z = \vec g$, see, e.g., \cite{DavisRajamanickamSidLakhdar2016, DuffErismanReid2017, GouldScottHu2007, LangerNeumuellerCism2018, Lui1992, Martinsson2020, PardoPaszynskiCollierAlvarezDalcinCalo2012}. Here, $K_{h_x} \in \C^{n \times n}$ is a sparse matrix coming from finite element/difference discretizations of a physical domain $\Omega \subset \R^d$, $d=1,2,3$, like the spatial mass or stiffness matrix \eqref{VF:Matritzen_Mass_Stiff}. Sparse direct solvers exploit the sparsity pattern of the system matrix $K_{h_x}$, and are based on a divide-and-conquer technique, which can be interpreted as procedure of subdividing the physical domain $\Omega$, which leads also to a subdivision of the degrees of freedom. Usually, the following steps have to be applied in such a method:
\begin{enumerate}
 \item Ordering Step, e.g., minimum degree or nested dissection methods, 
 \item Symbolic Factorization Step,
 \item Numerical Factorization Step,
 \item Solving Step.
\end{enumerate}
For structured grids, the complexity of these methods is summarized in Table~\ref{Vor:Tab:Loeser}.

\begin{table}[tbhp]
\begin{center}
\begin{tabular}{r|ccc} 
	\hline
 $d$ & Ordering and Factorization Steps & Solving Step & Memory  \\
		\hline   
    1 & $\mathcal O(n)$ & $\mathcal O(n)$ & $\mathcal O(n)$ \\
    2 & $\mathcal O(n^{3/2})$ & $\mathcal O(n \ln n)$ & $\mathcal O(n \ln n)$ \\
    3 & $\mathcal O(n^2)$ & $\mathcal O(n^{4/3})$ & $\mathcal O(n^{4/3})$ \\
    \hline
  \end{tabular}
    \caption{Summary of the complexity of sparse direct solver for sparse matrices coming from structured grids.} \label{Vor:Tab:Loeser}
\end{center}
\end{table}

There are several open-source software packages for sparse direct solvers. In this paper, we use only the sparse direct solver MUMPS 5.3.3 \cite{MUMPS2, MUMPS1}.

\section{Space-Time Solvers} \label{Sec:Loeser}

In this section, efficient solvers for the large-scale space-time system \eqref{VF:LGS} are developed. Our new solver is based on \cite[Section 3]{HoferLangerNeumuellerSchneckenleitner2019} and \cite[Section 4]{Tani2017}, where analogous results are derived for methods in isogeometric analysis. In greater detail, we state solvers for the global linear system 
\begin{equation} \label{Loeser:GlobalLGS}
  (A_{h_t}^{\mathcal H_T} \otimes M_{h_x} + M_{h_t}^{\mathcal H_T} \otimes A_{h_x}) \vec u = \vec{\widetilde F}^{\mathcal H_T}
\end{equation}
given in \eqref{VF:LGS} with the symmetric, positive definite matrices $M_{h_x} \in \R^{M_x \times M_x}$, $A_{h_x} \in \R^{M_x \times M_x}$, $A_{h_t}^{\mathcal H_T} \in \R^{N_t \times N_t}$ and the nonsymmetric, positive definite matrix $M_{h_t}^{\mathcal H_T} \in \R^{N_t \times N_t}$. Since \eqref{Loeser:GlobalLGS} is a (generalized) Sylvester equation, we can apply the Bartels-Stewart method \cite{BartelsStewart1972} with real- or complex-Schur decomposition and the Fast Diagonalization method \cite{Lynch1964} to solve \eqref{Loeser:GlobalLGS}, see also \cite{Gardiner1992, Simoncini2016}. In all three cases, the matrix pencil $(M_{h_t}^{\mathcal H_T}, A_{h_t}^{\mathcal H_T})$ is decomposed in the form
\begin{equation*}
    (A_{h_t}^{\mathcal H_T})^{-1} M_{h_t}^{\mathcal H_T} = X_t Z_t X_t^{-1}
\end{equation*}
with real, regular matrices $X_t, Z_t \in \R^{N_t \times N_t}$, where $Z_t$ is an upper (quasi-)triangular matrix, or complex, regular matrices $X_t, Z_t \in \C^{N_t \times N_t}$, where $Z_t$ is an upper triangular or diagonal matrix. Defining
\begin{equation*}
    Y_t := (A_{h_t}^{\mathcal H_T} X_t)^{-1}
\end{equation*}
gives the representations
\begin{equation*}
    A_{h_t}^{\mathcal H_T} = Y_t^{-1} X_t^{-1} \quad \text{ and } \quad  M_{h_t}^{\mathcal H_T} = \underbrace{ A_{h_t}^{\mathcal H_T}}_{= Y_t^{-1} X_t^{-1} } X_t Z_t X_t^{-1} = Y_t^{-1} Z_t X_t^{-1}.
\end{equation*}
Hence, the global linear system is equivalent to solving
\begin{equation*}
    (Y_t^{-1} \otimes I_{M_x}) (I_{N_t} \otimes M_{h_x} + Z_t \otimes A_{h_x}) (X_t^{-1} \otimes I_{M_x}) \vec u = \vec{\widetilde F}^{\mathcal H_T}
\end{equation*}
with the identity matrices $I_{N_t} \in \R^{N_t \times N_t}$ and $I_{M_x} \in \R^{M_x \times M_x}$. Thus, the solution of \eqref{Loeser:GlobalLGS} is given by
\begin{equation} \label{Loeser:LsgDarstellung}
    \vec u = (X_t \otimes I_{M_x}) (I_{N_t} \otimes M_{h_x} + Z_t \otimes A_{h_x})^{-1} (Y_t \otimes I_{M_x}) \vec{\widetilde F}^{\mathcal H_T}.
\end{equation}
The first step in \eqref{Loeser:LsgDarstellung} is the calculation of the vector
\begin{equation} \label{Loeser:Lsg_Schritt1}
 \vec g := (\vec g_1, \vec g_2, \ldots, \vec g_{N_t} )^\top := (Y_t \otimes I_{M_x}) \vec{\widetilde F}^{\mathcal H_T} =  \mathrm{vec} \left( \hat F  (A_{h_t}^{\mathcal H_T})^{-1} X_t^{-\top} \right) \in \C^{N_t\cdot M_x}
\end{equation}
with a matrix $\hat F \in \R^{M_x \times N_t}$ corresponding to the relation \eqref{Vor:VektorisierungVektor}, satisfying $\mathrm{vec}(\hat F) =  \vec{\widetilde F}^{\mathcal H_T},$ where $\vec g_i \in \C^{M_x}.$ The second step in \eqref{Loeser:LsgDarstellung} is to solve the linear system
\begin{equation*}
    (I_{N_t} \otimes M_{h_x} + Z_t \otimes A_{h_x})^{-1} \vec z = \vec g
\end{equation*}
for the vector
\begin{equation*}
    \vec z = (\vec z_1, \vec z_2, \ldots ,\vec z_{N_t})^\top \in \C^{N_t\cdot M_x},
\end{equation*}
where $\vec z_{\ell} \in \C^{M_x}$, $\ell=1,\dots,N_t$, which is analyzed in greater detail in the following subsections. The third step in \eqref{Loeser:LsgDarstellung} is the calculation of the desired unknown
\begin{equation} \label{Loeser:Lsg_Schritt3}
  \vec u = (X_t \otimes I_{M_x}) \vec z =  \mathrm{vec} \left( Z X_t^{\top} \right) \in \R^{N_t\cdot M_x}
\end{equation}
with a matrix $Z \in \C^{M_x \times N_t}$ corresponding to the relation \eqref{Vor:VektorisierungVektor}, satisfying $\mathrm{vec}(Z) =  \vec{z}.$

\subsection{Eigenvalues of the Matrix Pencil $(M_{h_t}^{\mathcal H_T}, A_{h_t}^{\mathcal H_T})$ }

In this subsection, we investigate the generalized eigenvalue problem
\begin{equation} \label{Loeser:EWP:EWP}
    M_{h_t}^{\mathcal H_T} \vec z = \lambda A_{h_t}^{\mathcal H_T} \vec z
\end{equation}
with eigenvalues $\lambda = \alpha + \mathrm{\iota} \beta \in \C$ and eigenvectors $\vec z = \vec x + \mathrm{\iota} \vec y \in \C^{N_t}$. As the matrix $M_{h_t}^{\mathcal H_T}$ is nonsymmetric, we have $\Im(\lambda) = \beta \neq 0$ in general, where the complex eigenvalues occur in conjugate pairs $\lambda = \alpha \pm \mathrm{\iota} \beta$. On the other hand, since the matrices $M_{h_t}^{\mathcal H_T}$ and $A_{h_t}^{\mathcal H_T}$ are positive definite,  and $A_{h_t}^{\mathcal H_T}$ is symmetric, it follows immediately from \cite[Lemma 3.2]{HoferLangerNeumuellerSchneckenleitner2019} that 
\begin{equation} \label{Loeser:EWP:positiverRT}
    \Re(\lambda) = \alpha > 0
\end{equation}
without any restriction on the mesh. Additionally, this property remains true for any conforming tensor-product ansatz space $\mathcal V_h$ in \eqref{VF:VF_HT_disk}, e.g., any polynomial degree $p \in \N$. In Subsection~\ref{Sec:Num:EW}, numerical examples, which investigate the eigenvalues $\lambda$, are given.

\subsection{Bartels-Stewart Method with Real-Schur Decomposition} \label{Sec:Loeser:BSr}

The aim of this subsection is to derive an algorithm on the basis of the Bartels-Stewart method with real-Schur decomposition \cite{BartelsStewart1972, Gardiner1992, Simoncini2016}. Therefore, a real-Schur decomposition of the matrix pencil $(M_{h_t}^{\mathcal H_T}, A_{h_t}^{\mathcal H_T})$ is used in the form
\begin{equation} \label{Loeser:BSr:Schurzerlegung}
    (A_{h_t}^{\mathcal H_T})^{-1} M_{h_t}^{\mathcal H_T} = Q_t R_t Q_t^\top
\end{equation}
with the orthogonal matrix $Q_t := X_t \in \R^{N_t \times N_t}$ and the upper quasi-triangular matrix $R_t := Z_t \in \R^{N_t \times N_t},$ where the diagonals of $R_t$ have $2\times2$ and $1\times1$ blocks, corresponding to complex and real eigenvalues of the matrix $(A_{h_t}^{\mathcal H_T})^{-1} M_{h_t}^{\mathcal H_T}$. In greater detail, let $\lambda_1,\dots,\lambda_{N_t} \in \C$ be the eigenvalues of the matrix $ (A_{h_t}^{\mathcal H_T})^{-1} M_{h_t}^{\mathcal H_T} \in \R^{N_t\times N_t}$. To each real eigenvalue $\lambda \in \R$, we can relate a $1\times1$ block given as $\begin{pmatrix} \lambda \end{pmatrix} \in \R^{1\times 1}$. The complex eigenvalues occur in conjugate pairs. Thus, each conjugate pair $\alpha \pm \mathrm{\iota} \beta$ corresponds to  a $2\times2$ block
\begin{equation*}
 \begin{pmatrix}
  \alpha & b_1 \\
  b_2    & \alpha
 \end{pmatrix} \in \R^{2\times2},
\end{equation*}
satisfying $\abs{\beta} = \sqrt{\abs{b_1 b_2}}>0$ with $b_1$ and $b_2$ having different signs. With the real-Schur decomposition \eqref{Loeser:BSr:Schurzerlegung} and \eqref{Loeser:LsgDarstellung}, the solution of \eqref{Loeser:GlobalLGS} can be represented in the form
\begin{equation*}
    \vec u = (Q_t \otimes I_{M_x}) (I_{N_t} \otimes M_{h_x} + R_t \otimes A_{h_x})^{-1} (Y_t \otimes I_{M_x}) \vec{\widetilde F}^{\mathcal H_T},
\end{equation*}
where $Y_t =  Q_t^\top (A_{h_t}^{\mathcal H_T})^{-1}$. The applications of the transformation matrices $Y_t \otimes I_{M_x} = Q_t^\top (A_{h_t}^{\mathcal H_T})^{-1} \otimes I_{M_x}$ and $Q_t \otimes I_{M_x}$ are given by \eqref{Loeser:Lsg_Schritt1} and \eqref{Loeser:Lsg_Schritt3}. Hence, it remains to solve
\begin{equation} \label{Loeser:BSr:LGS}
    (I_{N_t} \otimes M_{h_x} + R_t \otimes A_{h_x}) \vec z = \vec g := (Y_t \otimes I_{M_x}) \vec{\widetilde F}^{\mathcal H_T}
\end{equation}
with the unknown
\begin{equation*}
    \vec z = (\vec z_1, \vec z_2, \ldots, \vec z_{N_t} )^\top \in \R^{N_t\cdot M_x},
\end{equation*}
where $\vec z_{\ell} \in \R^{M_x}$, $\ell=1,\dots,N_t$. In addition, let the vector of the right-hand side
\begin{equation*}
    \vec g = (\vec g_1, \vec g_2, \ldots, \vec g_{N_t})^\top \in \R^{N_t\cdot M_x}
\end{equation*}
be decomposed, where $\vec g_{\ell} \in \R^{M_x}$, $\ell=1,\dots,N_t$. Since the global linear system \eqref{Loeser:BSr:LGS} has a special triangular structure, this system can be solved by a backward substitution technique, which is described in the following in more detail. Therefore, let $k \in \{1,\dots,N_t-1\}$ be such that $\vec z_{k+1}, \dots, \vec z_{N_t}$ are already computed, or let $k=N_t$, where we set $\sum_{j=N_t+1}^{N_t} (\cdot) := 0.$ Then, two cases occur, as the diagonals of $R_t$ have $2\times2$ and $1\times1$ blocks:
\begin{enumerate}
 \item In the case of a $1\times1$ block of $R_t$, i.e. $k=1$ or $R_t[k-1,k] = 0$, the linear system
    \begin{equation} \label{Loeser:BSr:LGS_Einserblock}
        (M_{h_x} + R_t[k,k] A_{h_x}) \vec z_k = \vec g_k - \sum_{j=k+1}^{N_t} R_t[k,j] A_{h_x} \vec z_j
    \end{equation}
    has to be solved for $\vec z_k.$
 \item In the case of a $2\times2$ block of $R_t$, i.e. $k>1$ and $R_t[k-1,k] \neq 0$, the linear system
    \begin{multline} \label{Loeser:BSr:LGS_Zweierblock}
        \begin{pmatrix}
            M_{h_x} + \alpha A_{h_x} &  b_1 A_{h_x} \\
            b_2 A_{h_x}               & M_{h_x} + \alpha A_{h_x}
        \end{pmatrix}
        \begin{pmatrix}
            \vec z_{k-1} \\
            \vec z_k
        \end{pmatrix}
        \\
        =
        \begin{pmatrix}
         \vec g_{k-1} - \sum_{j=k+1}^{N_t} R_t[k-1,j] A_{h_x} \vec z_j \\
         \vec g_k - \sum_{j=k+1}^{N_t} R_t[k,j] A_{h_x} \vec z_j
        \end{pmatrix}  
    \end{multline}
     with
     \begin{equation*}
        \alpha =  R_t[k-1,k-1] =  R_t[k,k]>0, \quad b_1 = R_t[k-1,k] \neq 0, \quad b_2 = R_t[k,k-1] \neq 0
     \end{equation*}
      has to be solved for $\vec z_{k-1}$ and $\vec z_k$.
\end{enumerate}
The system matrix of the linear system \eqref{Loeser:BSr:LGS_Einserblock} is symmetric and positive definite, since $R_t[k,k]>0$ is a real eigenvalue of the matrix $ (A_{h_t}^{\mathcal H_T})^{-1} M_{h_t}^{\mathcal H_T}$. The linear system \eqref{Loeser:BSr:LGS_Zweierblock} is equivalent to the linear system
\begin{multline} \label{Loeser:BSr:LGS_Zweierblock_symmetrisch}
    \begin{pmatrix}
        \abs{b_2} \left( M_{h_x} + \alpha A_{h_x} \right) &  -b_1 \abs{b_2} A_{h_x} \\
        \abs{b_1} b_2 A_{h_x}                             & -\abs{b_1} \left( M_{h_x} + \alpha A_{h_x} \right)
    \end{pmatrix}
    \begin{pmatrix}
        \vec z_{k-1} \\
        -\vec z_k
    \end{pmatrix}
    \\
    =
    \begin{pmatrix}
        \abs{b_2} \left( \vec g_{k-1} - \sum_{j=k+1}^{N_t} R_t[k-1,j] A_{h_x} \vec z_j \right)\\
        \abs{b_1} \left( \vec g_k - \sum_{j=k+1}^{N_t} R_t[k,j] A_{h_x} \vec z_j \right)
    \end{pmatrix},  
\end{multline}
where the system matrix is symmetric, but indefinite due to the property $b_1 \abs{b_2} = - \abs{b_1} b_2$. Note that the linear systems \eqref{Loeser:BSr:LGS_Zweierblock}, \eqref{Loeser:BSr:LGS_Zweierblock_symmetrisch} are uniquely solvable, since multiplying the second equation in \eqref{Loeser:BSr:LGS_Zweierblock_symmetrisch} by $-1$ leads to a nonsymmetric, but positive definite system matrix. The spatial linear systems \eqref{Loeser:BSr:LGS_Einserblock} and \eqref{Loeser:BSr:LGS_Zweierblock_symmetrisch} can be solved by (preconditioned) iterative solvers or by direct solvers. In this work, we consider sparse direct solvers only. The resulting algorithm of the Bartels-Stewart method with real-Schur decomposition is summarized in Algorithm~\ref{Loeser:BSr:Alg}.

\begin{algorithm} 
\caption{Bartels-Stewart method with real-Schur decomposition with output $\vec u$.} \label{Loeser:BSr:Alg}
\begin{algorithmic}
    \item[1:] Compute the real-Schur decomposition $(A_{h_t}^{\mathcal H_T})^{-1} M_{h_t}^{\mathcal H_T} = Q_t R_t Q_t^\top$ in \eqref{Loeser:BSr:Schurzerlegung}.
    \item[2:] Solve $\vec g = (\vec g_1, \vec g_2, \dots, \vec g_{N_t})^\top := (Y_t \otimes I_{M_x}) \vec{\widetilde F}^{\mathcal H_T} =  \mathrm{vec} \left( \hat F  (A_{h_t}^{\mathcal H_T})^{-1} Q_t \right)$  in \eqref{Loeser:Lsg_Schritt1}.
    \item[3:] Set $k=N_t$ and compute $\vec z = (\vec z_1, \vec z_2, \dots, \vec z_{N_t})^\top$ sequentially
        by
        \WHILE{$k>0$}
            \IF{$k=1$ \OR $R_t[k-1,k] = 0$}
                \STATE{Solve
                \begin{equation*}
                    (M_{h_x} + R_t[k,k] A_{h_x}) \vec z_k = \vec g_k - \sum_{j=k+1}^{N_t} R_t[k,j] A_{h_x} \vec z_j
                \end{equation*}
                for $\vec z_k$ in \eqref{Loeser:BSr:LGS_Einserblock}.}
                Set $k=k-1.$
            \ELSE
            \STATE{ Set
                \begin{equation*}
                    \alpha = R_t[k-1,k-1], \quad b_1 = R_t[k-1,k], \quad b_2 = R_t[k,k-1]
                \end{equation*}
                and solve
                \begin{multline*}
                    \begin{pmatrix}
                        \abs{b_2} \left( M_{h_x} + \alpha A_{h_x} \right) &  -b_1 \abs{b_2} A_{h_x} \\
                        \abs{b_1} b_2 A_{h_x}                             & -\abs{b_1} \left( M_{h_x} + \alpha A_{h_x} \right)
                    \end{pmatrix}
                    \begin{pmatrix}
                        \vec z_{k-1} \\
                        -\vec z_k
                    \end{pmatrix}
                    \\
                    =
                    \begin{pmatrix}
                        \abs{b_2} \left( \vec g_{k-1} - \sum_{j=k+1}^{N_t} R_t[k-1,j] A_{h_x} \vec z_j \right)\\
                        \abs{b_1} \left( \vec g_k - \sum_{j=k+1}^{N_t} R_t[k,j] A_{h_x} \vec z_j \right)
                    \end{pmatrix}
                \end{multline*}
                for  $(\vec z_{k-1}, -\vec z_k)^\top$ in \eqref{Loeser:BSr:LGS_Zweierblock_symmetrisch}. Set $k=k-2.$
            }
            \ENDIF
        \ENDWHILE
    \item[4:] Compute the matrix-vector product $\vec u = (Q_t \otimes I_{M_x}) \vec z =  \mathrm{vec} \left( Z Q_t^\top \right)$ in \eqref{Loeser:Lsg_Schritt3}.
\end{algorithmic}
\end{algorithm}

Numerical examples, which investigate the Bartels-Stewart method with real-Schur decomposition, are given in Subsection~\ref{Sec:Num:BSr}.

\subsubsection{Computational Cost and Memory Requirement}

The computational cost of step 1 in Algorithm~\ref{Loeser:BSr:Alg}, i.e. the real-Schur decomposition in \eqref{Loeser:BSr:Schurzerlegung}, is $\mathcal O( N_t^3)$, whereas the memory demand is $\mathcal O(N_t^2).$
To perform step 2 in Algorithm~\ref{Loeser:BSr:Alg} for calculating the vector $\vec g$ in \eqref{Loeser:Lsg_Schritt1}, $M_x$ linear systems of the size $N_t$ have to be solved for the same system matrix $A_{h_t}^{\mathcal H_T}$, and a matrix multiplication with $Q_t$ has to be applied. Using a Cholesky factorization of $A_{h_t}^{\mathcal H_T}$ of costs $\mathcal O(N_t^3)$ yields total computational costs of $\mathcal O(M_x N_t^2 + N_t^3)$ and a memory demand of $\mathcal O(N_t^2 + M_x N_t)$ for step 2 in Algorithm~\ref{Loeser:BSr:Alg}. Also step 4 in Algorithm~\ref{Loeser:BSr:Alg} requires computational costs of $\mathcal O(M_x N_t^2)$ and a memory consumption of $\mathcal O(N_t^2 + M_x N_t)$. The most expensive part of Algorithm~\ref{Loeser:BSr:Alg} is step 3, i.e. solving $N_t$  spatial linear systems of the form \eqref{Loeser:BSr:LGS_Einserblock} or \eqref{Loeser:BSr:LGS_Zweierblock_symmetrisch}. Assume that solving a spatial linear systems of the form \eqref{Loeser:BSr:LGS_Einserblock} or \eqref{Loeser:BSr:LGS_Zweierblock_symmetrisch} requires $\mathcal O(C_C(M_x))$ operations and $\mathcal O(C_S(M_x))$ storage with the cost function $C_C(\cdot)$ and the storage function $C_S(\cdot)$ defined by the spatial solver for the corresponding linear systems. Then, step 3 costs
$\mathcal O(C_C(M_x) \cdot N_t)$ for computations and $\mathcal O(C_S(M_x) + M_x N_t)$ for storage, where the calculation of the right-hand sides in \eqref{Loeser:BSr:LGS_Einserblock} or \eqref{Loeser:BSr:LGS_Zweierblock_symmetrisch} is of costs of lower order due to $A_{h_x}$ is sparse.
Hence, the overall computational cost and memory consumption of Algorithm~\ref{Loeser:BSr:Alg} are
\begin{equation} \label{Loeser:BSr:Komplexitaet}
 \mathcal O(N_t^3 +  M_x N_t^2 + C_C(M_x) \cdot N_t) \quad \text{ and } \quad \mathcal O(N_t^2 + M_x N_t + C_S(M_x)).
\end{equation}
Note that the calculations corresponding to the term $M_x N_t^2$ are few matrix multiplications, which are parallelizable and can be written as highly efficient BLAS-3 operations.
For the case of a uniform refinement strategy in temporal and spatial direction, i.e. $N_t$ doubles and $M_x$ grows by a factor $\mathcal O(2^d)$ in each refinement step, the number of the degrees of freedom $\mathrm{dof} = N_t \cdot M_x$ increases by a factor $\mathcal O(2^{d+1}).$ Hence, we have $N_t \sim M_x^{1/d} \sim \mathrm{dof}^{1/(d+1)}$, which results in the complexity of Algorithm~\ref{Loeser:BSr:Alg}, given in Table~\ref{Loeser:BSr:Tab:Komplexitaet}, when a sparse direct solver of Subsection~\ref{Sec:Vor:SparseLoeser} is applied for the spatial problems of step 3 in Algorithm~\ref{Loeser:BSr:Alg} in the case of structured grids. Note that the sparsity patterns of the system matrices in step 3 of Algorithm~\ref{Loeser:BSr:Alg} remain the same for $k=1,\dots,N_t$, which can be exploited by the sparse direct solver of Subsection~\ref{Sec:Vor:SparseLoeser}, i.e. it is sufficient to perform the ordering and symbolic factorization steps only once.

\begin{table}[tbhp]
\begin{center}
\begin{tabular}{r|cc} 
	\hline
 $d$ & Computations & Memory  \\
		\hline   
    1 & $\mathcal O(\mathrm{dof}^{3/2})$ & $\mathcal O(\mathrm{dof})$ \\
    2 & $\mathcal O(\mathrm{dof}^{4/3})$ & $\mathcal O(\mathrm{dof})$ \\
    3 & $\mathcal O(\mathrm{dof}^{7/4})$ & $\mathcal O(\mathrm{dof})$ \\
    \hline
  \end{tabular}
    \caption{Summary of the complexity of the Bartels-Stewart methods (Algorithm~\ref{Loeser:BSr:Alg}, Algorithm~\ref{Loeser:BSc:Alg}) or the Fast Diagonalization method (Algorithm~\ref{Loeser:FD:Alg}), using a sparse direct solver for spatial structured grids given in Subsection~\ref{Sec:Vor:SparseLoeser}.} \label{Loeser:BSr:Tab:Komplexitaet}
\end{center}
\end{table}

\subsection{Bartels-Stewart Method with Complex-Schur Decomposition} \label{Sec:Loeser:BSc}

In this subsection, an algorithm, using the Bartels-Stewart method with complex-Schur decomposition \cite{BartelsStewart1972, Gardiner1992, Simoncini2016}, is derived. Therefore, a complex-Schur decomposition of the matrix pencil $(M_{h_t}^{\mathcal H_T}, A_{h_t}^{\mathcal H_T})$ is used in the form
\begin{equation} \label{Loeser:BSc:Schurzerlegung}
    (A_{h_t}^{\mathcal H_T})^{-1} M_{h_t}^{\mathcal H_T} = W_t S_t W_t^*
\end{equation}
with the unitary matrix $W_t := X_t \in \C^{N_t \times N_t}$ and the upper triangular matrix $S_t := Z_t \in \C^{N_t \times N_t},$ where the generalized eigenvalues $\lambda_1,\dots,\lambda_{N_t} \in \C$ of the matrix pencil $(M_{h_t}^{\mathcal H_T}, A_{h_t}^{\mathcal H_T})$ are on the diagonal of $S_t$, i.e. $S_t[k,k] = \lambda_k$ for $k=1,\dots,N_t$. Note that the real parts of the complex eigenvalues fulfill $\Re(\lambda_\ell) > 0$ for all $\ell=1,\dots,N_t,$ see \eqref{Loeser:EWP:positiverRT}. With the complex-Schur decomposition \eqref{Loeser:BSc:Schurzerlegung} and \eqref{Loeser:LsgDarstellung}, the solution of \eqref{Loeser:GlobalLGS} is given by
\begin{equation*}
    \vec u = (W_t \otimes I_{M_x}) (I_{N_t} \otimes M_{h_x} + S_t \otimes A_{h_x})^{-1} (Y_t \otimes I_{M_x}) \vec{\widetilde F}^{\mathcal H_T},
\end{equation*}
where $Y_t =  W_t^* (A_{h_t}^{\mathcal H_T})^{-1}$. The applications of the transformation matrices $Y_t \otimes I_{M_x} = W_t^* (A_{h_t}^{\mathcal H_T})^{-1} \otimes I_{M_x}$ and $W_t \otimes I_{M_x}$ are given by \eqref{Loeser:Lsg_Schritt1} and \eqref{Loeser:Lsg_Schritt3}. Hence, it remains to solve
\begin{equation} \label{Loeser:BSc:LGS}
    (I_{N_t} \otimes M_{h_x} + S_t \otimes A_{h_x}) \vec z = \vec g := (Y_t \otimes I_{M_x}) \vec{\widetilde F}^{\mathcal H_T}
\end{equation}
with the unknown
\begin{equation*}
    \vec z = (\vec z_1, \vec z_2, \ldots, \vec z_{N_t} )^\top \in \C^{N_t\cdot M_x},
\end{equation*}
where $\vec z_{\ell} \in \C^{M_x}$, $\ell=1,\dots,N_t$. In addition, let the vector of the right-hand side
\begin{equation*}
    \vec g = (\vec g_1, \vec g_2, \ldots, \vec g_{N_t})^\top \in \C^{N_t\cdot M_x}
\end{equation*}
be decomposed, where $\vec g_{\ell} \in \C^{M_x}$, $\ell=1,\dots,N_t$. Since the global linear system \eqref{Loeser:BSc:LGS} has a block triangular structure, this system can be solved by a backward substitution technique, which is described in the following in more detail.
For $k=N_t,N_t-1,\dots,1$, the linear system
\begin{equation} \label{Loeser:BSc:LGS_x}
    (M_{h_x} + \underbrace{S_t[k,k]}_{=\lambda_k} A_{h_x}) \vec z_k = \vec g_k - \sum_{j=k+1}^{N_t} S_t[k,j] A_{h_x} \vec z_j
\end{equation}
has to be solved for $\vec z_k$, where we set $\sum_{j=N_t+1}^{N_t} (\cdot) := 0.$ Note that the system matrix of the linear system \eqref{Loeser:BSc:LGS_x} is symmetric, but not Hermitian for $\Im(S_t[k,k]) = \Im(\lambda_k) \neq 0$. With separating the real and the imaginary parts $\lambda_k = \alpha_k + \iota \beta_k$, the complex linear system \eqref{Loeser:BSc:LGS_x} is equivalent to a real linear system of doubled size with a system matrix
\begin{equation*}
  \begin{pmatrix}
    M_{h_x} + \alpha_k A_{h_x} & -\beta_k A_{h_x} \\
    \beta_k A_{h_x}            &  M_{h_x} + \alpha_k A_{h_x}
  \end{pmatrix},
\end{equation*}
which is nonsymmetric, but positive definite due to $\alpha_k >0$, see \eqref{Loeser:EWP:positiverRT}. Hence, the spatial linear systems \eqref{Loeser:BSc:LGS_x} are uniquely solvable, which can be solved by (preconditioned) iterative solvers or by direct solvers. 
In this work, we consider sparse direct solvers only. The resulting algorithm of the Bartels-Stewart method with complex-Schur decomposition is summarized in Algorithm~\ref{Loeser:BSc:Alg}, where $\overline{W_t}$ is the element-by-element conjugation of $W_t.$

\begin{algorithm} 
\caption{Bartels-Stewart method with complex-Schur decomposition with output $\vec u$.} \label{Loeser:BSc:Alg}
\begin{algorithmic}
    \item[1:] Compute the complex-Schur decomposition $(A_{h_t}^{\mathcal H_T})^{-1} M_{h_t}^{\mathcal H_T} = W_t S_t W_t^*$ in \eqref{Loeser:BSc:Schurzerlegung}.
    \item[2:] Solve $\vec g = (\vec g_1, \vec g_2, \dots, \vec g_{N_t})^\top := (Y_t \otimes I_{M_x}) \vec{\widetilde F}^{\mathcal H_T} =  \mathrm{vec} \left( \hat F  (A_{h_t}^{\mathcal H_T})^{-1} \overline{W_t} \right)$  in \eqref{Loeser:Lsg_Schritt1}.
    \item[3:] Compute $\vec z = (\vec z_1, \vec z_2, \dots, \vec z_{N_t})^\top$ sequentially
        by
        \FOR{$k=N_t, N_t-1,\ldots,1$}
            \STATE{Solve
            \begin{equation*}
                (M_{h_x} + S_t[k,k] A_{h_x}) \vec z_k = \vec g_k - \sum_{j=k+1}^{N_t} S_t[k,j] A_{h_x} \vec z_j
            \end{equation*}
            for $\vec z_k$ in \eqref{Loeser:BSc:LGS_x}.}
        \ENDFOR
    \item[4:] Compute the matrix-vector product $\vec u = (W_t \otimes I_{M_x}) \vec z =  \mathrm{vec} \left( Z W_t^\top \right)$ in \eqref{Loeser:Lsg_Schritt3}.
\end{algorithmic}
\end{algorithm}

Numerical examples, which investigate the Bartels-Stewart method with complex-Schur decomposition, are given in Subsection~\ref{Sec:Num:BSc}.

\subsubsection{Computational Cost and Memory Requirement}

The computational cost and memory requirement of Algorithm~\ref{Loeser:BSc:Alg} can be analyzed in the same way as for the Bartels-Stewart method with real-Schur decomposition (Algorithm~\ref{Loeser:BSr:Alg}). Hence, the overall computational cost and memory consumption of Algorithm~\ref{Loeser:BSc:Alg} are
\begin{equation*}
 \mathcal O(N_t^3 +  M_x N_t^2 + C_C(M_x) \cdot N_t) \quad \text{ and } \quad \mathcal O(N_t^2 + M_x N_t + C_S(M_x)),
\end{equation*}
which are of the same order as for the Bartels-Stewart method with real-Schur decomposition, see \eqref{Loeser:BSr:Komplexitaet}. Note that the calculations corresponding to the term $M_x N_t^2$ are few matrix multiplications, which are parallelizable and can be written as highly efficient BLAS-3 operations.
For the case of a uniform refinement strategy in temporal and spatial direction for spatial structured grids, the complexity of Algorithm~\ref{Loeser:BSc:Alg} is again given in Table~\ref{Loeser:BSr:Tab:Komplexitaet}. Note that the sparsity patterns of the system matrices in step 3 of Algorithm~\ref{Loeser:BSc:Alg} remain the same for $k=1,\dots,N_t$, which can be exploited by the sparse direct solver of Subsection~\ref{Sec:Vor:SparseLoeser}, i.e. it is sufficient to perform the ordering and symbolic factorization steps only once.

\subsection{Fast Diagonalization Method} \label{Sec:Loeser:FD}

This subsection deals with the development of an algorithm that is based on the Fast Diagonalization method \cite{Lynch1964, Simoncini2016}. Therefore, an eigenvalue decomposition of the matrix pencil $(M_{h_t}^{\mathcal H_T}, A_{h_t}^{\mathcal H_T})$ is used in the form
\begin{equation} \label{Loeser:FD:Zerlegung}
    (A_{h_t}^{\mathcal H_T})^{-1} M_{h_t}^{\mathcal H_T} = X_t D_t X_t^{-1}
\end{equation}
with the complex matrix $X_t \in \C^{N_t \times N_t}$ of generalized eigenvectors and the complex diagonal matrix
\begin{equation*}
    D_t := Z_t := \mathrm{diag}(\lambda_{1},\dots,\lambda_{N_t}) \in \C^{N_t \times N_t}
\end{equation*}
with the complex generalized eigenvalues $\lambda_\ell \in \C$, $\ell=1,\dots,N_t.$ The real parts of the complex eigenvalues fulfill $\Re(\lambda_\ell) > 0$ for all $\ell=1,\dots,N_t,$ see \eqref{Loeser:EWP:positiverRT}. Since the matrix $M_{h_t}^{\mathcal H_T}$ is nonsymmetric, the matrix $X_t$ of generalized eigenvectors is not unitary and so, its condition number is not 1. As the condition number of $X_t$ may be large, numerical instabilities may occur by applying the inverse of $X_t$, which may be damped by an additional singular decomposition of $X_t$. Hence, we apply the singular value decomposition
\begin{equation} \label{Loeser:FD:SVD}
 X_t = U_t \Sigma_t V_t^*
\end{equation}
with unitary matrices $U_t, V_t \in \C^{N_t \times N_t}$ and the diagonal matrix $\Sigma_t \in \R^{N_t \times N_t}$.  With the diagonalization \eqref{Loeser:FD:Zerlegung} and \eqref{Loeser:LsgDarstellung}, the solution of \eqref{Loeser:GlobalLGS} is given by
\begin{equation*}
    \vec u = (X_t \otimes I_{M_x}) (I_{N_t} \otimes M_{h_x} + D_t \otimes A_{h_x})^{-1} (Y_t \otimes I_{M_x}) \vec{\widetilde F}^{\mathcal H_T}.
\end{equation*}
With the singular value decomposition \eqref{Loeser:FD:SVD}, the representation
\begin{equation*}
 Y_t = X_t^{-1} (A_{h_t}^{\mathcal H_T})^{-1} = V_t \Sigma_t^{-1} U_t^* (A_{h_t}^{\mathcal H_T})^{-1}
\end{equation*}
gives the transformation matrices $Y_t \otimes I_{M_x} = V_t \Sigma_t^{-1} U_t^* (A_{h_t}^{\mathcal H_T})^{-1} \otimes I_{M_x}$ and $X_t \otimes I_{M_x} = U_t \Sigma_t V_t^* \otimes I_{M_x}$ for the calculations in \eqref{Loeser:Lsg_Schritt1} and \eqref{Loeser:Lsg_Schritt3}. Hence, it remains to solve $N_t$ spatial problems with the complex system matrix
\begin{equation} \label{Loeser:FD:LGS_unabhaengig}
    M_{h_x} + \lambda_\ell A_{h_x} \in \C^{M_x \times M_x}
\end{equation}
for $\ell=1,\dots,N_t$, which can be done independently, i.e. a parallelization in the time direction is possible.  The system matrices \eqref{Loeser:FD:LGS_unabhaengig} are the same as in \eqref{Loeser:BSc:LGS_x} for the Bartels-Stewart method with complex-Schur decomposition, i.e. they are regular and symmetric, but not Hermitian for $\Im(\lambda_k) \neq 0$. The spatial linear systems with the system matrix \eqref{Loeser:FD:LGS_unabhaengig} can be solved by (preconditioned) iterative solvers or by direct solvers. In this work, we consider sparse direct solvers only. The resulting algorithm of the Fast Diagonalization method is summarized in Algorithm~\ref{Loeser:FD:Alg}, where $\overline{U_t}, \overline{V_t}$ are the element-by-element conjugations of $U_t, V_t.$

\begin{algorithm} 
\caption{Fast Diagonalization method with output $\vec u$.} \label{Loeser:FD:Alg}
\begin{algorithmic}
    \item[1a:] Compute the eigenvalue decomposition $(A_{h_t}^{\mathcal H_T})^{-1} M_{h_t}^{\mathcal H_T} = X_t D_t X_t^{-1}$ in \eqref{Loeser:FD:Zerlegung}.
    \item[1b:] Compute the singular value decomposition $X_t = U_t \Sigma_t V_t^*$ in \eqref{Loeser:FD:SVD}.
    \item[2:] Solve $\vec g = (\vec g_1, \vec g_2, \dots, \vec g_{N_t})^\top := (Y_t \otimes I_{M_x}) \vec{\widetilde F}^{\mathcal H_T} =  \mathrm{vec} \left( \hat F  (A_{h_t}^{\mathcal H_T})^{-1} \overline{U_t} \Sigma_t^{-1} V_t^\top \right)$  in \eqref{Loeser:Lsg_Schritt1}.
    \item[3:] Compute $\vec z = (\vec z_1, \vec z_2, \dots, \vec z_{N_t})^\top$ in parallel by
        \FOR{$k=1,\dots,N_t$}
                \STATE{Solve
                \begin{equation*}
                    (M_{h_x} + \lambda_k A_{h_x}) \vec z_k = \vec g_k
                \end{equation*}
                for $\vec z_k$.}
        \ENDFOR
    \item[4:] Compute the matrix-vector product $\vec u = (U_t \Sigma_t V_t^* \otimes I_{M_x}) \vec z = \mathrm{vec} \left( Z \overline{V_t} \Sigma_t U_t^\top \right)$ in \eqref{Loeser:Lsg_Schritt3}.
\end{algorithmic}
\end{algorithm}

Numerical examples, which investigate the Fast Diagonalization method, are given in Subsection~\ref{Sec:Num:FD}.

\subsubsection{Computational Cost and Memory Requirement}

The computational cost of step 1a and step 1b in Algorithm~\ref{Loeser:FD:Alg}, i.e. the eigenvalue decomposition in \eqref{Loeser:FD:Zerlegung} and the singular value decomposition \eqref{Loeser:FD:SVD}, is $\mathcal O( N_t^3)$, whereas the memory demand is $\mathcal O(N_t^2).$ To perform step 2 in Algorithm~\ref{Loeser:FD:Alg} for calculating the vector $\vec g$ in \eqref{Loeser:Lsg_Schritt1}, $M_x$ linear systems of the size $N_t$ have to be solved for the same system matrix $A_{h_t}^{\mathcal H_T}$, and matrix multiplications with $\overline{U_t} \Sigma_t^{-1} V_t^\top$ have to be applied. Using a Cholesky factorization of $A_{h_t}^{\mathcal H_T}$ of costs $\mathcal O(N_t^3)$ yields total computational costs of $\mathcal O(M_x N_t^2 + N_t^3)$ and a memory demand of $\mathcal O(N_t^2 + M_x N_t)$ for step 2 in Algorithm~\ref{Loeser:FD:Alg}. Also step 4 in Algorithm~\ref{Loeser:FD:Alg} requires computational costs of $\mathcal O(M_x N_t^2)$ and a memory consumption of $\mathcal O(N_t^2 + M_x N_t)$. The most expensive part of Algorithm~\ref{Loeser:FD:Alg} is step 3, i.e. solving $N_t$ spatial linear systems with the system matrix \eqref{Loeser:FD:LGS_unabhaengig}, which can be done in parallel. Assume that solving a spatial linear systems with the system matrix \eqref{Loeser:FD:LGS_unabhaengig} requires $\mathcal O(C_C(M_x))$ operations and $\mathcal O(C_S(M_x))$ storage with some cost function $C_C(\cdot)$ and some storage function $C_S(\cdot)$. Then, step 3 costs $\mathcal O(C_C(M_x) \cdot N_t)$ for computations and $\mathcal O(C_S(M_x) + M_x N_t)$ for storage. Hence, the overall computational cost and memory consumption of Algorithm~\ref{Loeser:FD:Alg} are
\begin{equation*}
 \mathcal O(N_t^3 +  M_x N_t^2 + C_C(M_x) \cdot N_t) \quad \text{ and } \quad \mathcal O(N_t^2 + M_x N_t + C_S(M_x)),
\end{equation*}
which are of the same order as for the Bartels-Stewart method with real-Schur decomposition, see \eqref{Loeser:BSr:Komplexitaet}. Note that the calculations corresponding to the term $M_x N_t^2$ are few matrix multiplications, which are parallelizable and can be written as highly efficient BLAS-3 operations.
For the case of a uniform refinement strategy in temporal and spatial direction for spatial structured grids, the complexity of Algorithm~\ref{Loeser:FD:Alg} is again given in Table~\ref{Loeser:BSr:Tab:Komplexitaet}. Note that the sparsity patterns of the system matrices in step 3 of Algorithm~\ref{Loeser:FD:Alg} remain the same for $k=1,\dots,N_t$, which can be exploited by the sparse direct solver of Subsection~\ref{Sec:Vor:SparseLoeser}, i.e. it is sufficient to perform the ordering and symbolic factorization steps only once.

\section{Numerical Examples} \label{Sec:Num}

\begin{figure}[tbhp]
\begin{center}
    \includegraphics[scale=0.47]{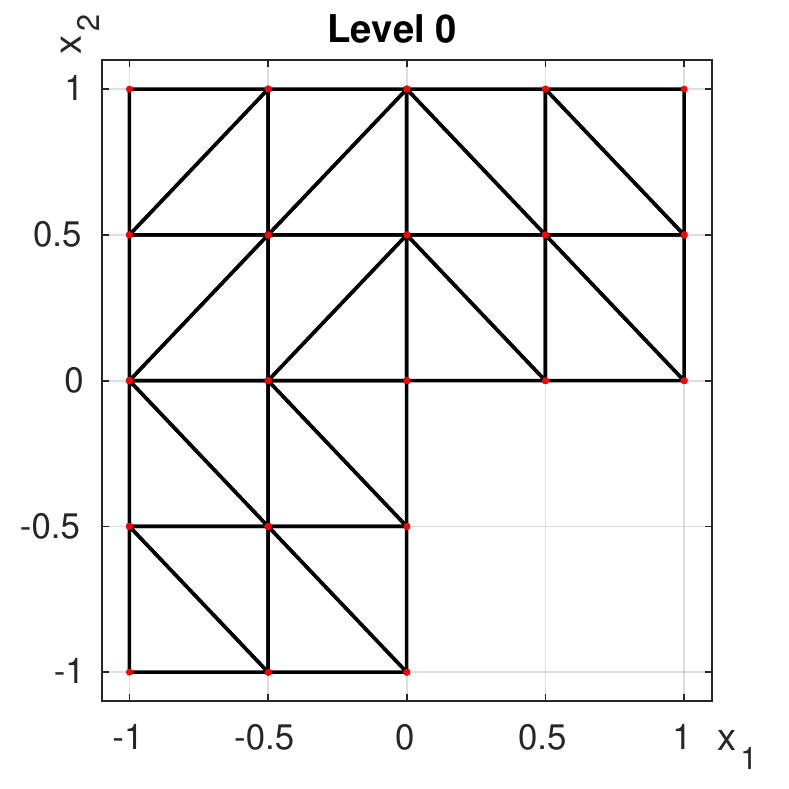}
    \includegraphics[scale=0.47]{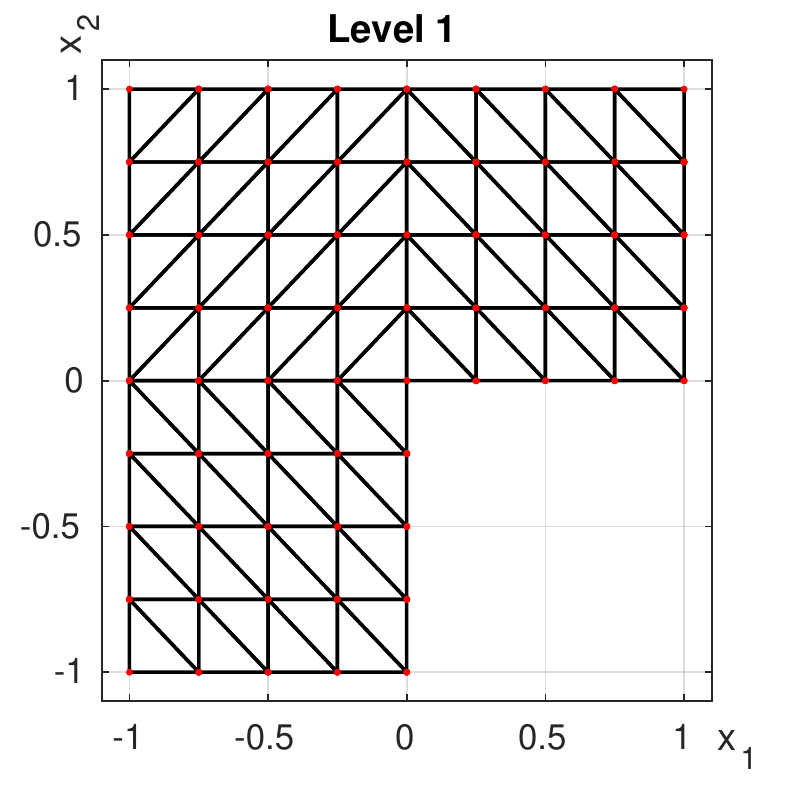}
    \includegraphics[scale=0.47]{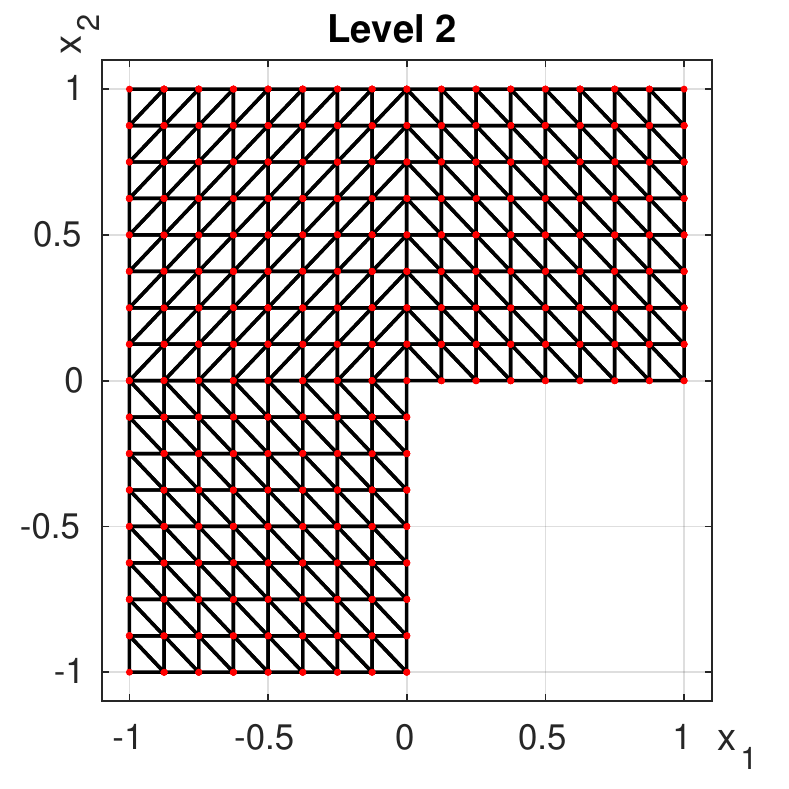}
  \caption{Uniform refinement strategy: Starting mesh, the meshes after one and two uniform refinement steps.}
  \label{Num:Fig:NetzeDreieckGlm}
\end{center}
\end{figure}

In this section, numerical examples for the generalized eigenvalue problem \eqref{Loeser:EWP:EWP} and for the Galerkin finite element method \eqref{VF:VF_HT_disk_Tensor_gestoert} using the Bartels-Stewart methods (Algorithm~\ref{Loeser:BSr:Alg} and Algorithm~\ref{Loeser:BSc:Alg}) and the Fast Diagonalization method (Algorithm~\ref{Loeser:FD:Alg}) are given. As numerical example, we consider the parabolic initial-boundary value problem \eqref{Einf:WaermePDG} in the two-dimensional spatial L-shaped domain
\begin{equation} \label{Num:Lshape}
  \Omega := (-1,1)^2 \setminus \left( [0,1] \times [-1,0] \right) \subset \R^2
\end{equation}
and with the terminal time $T=\frac 1 2$. We use the manufactured solution
\begin{equation} \label{Num:ExakteLsg}
  u(x_1,x_2,t) = \frac{5}{2 \pi  t} \mathrm e^{\frac{-\left({x_1}-\frac{1}{4}\right)^2-\left({x_2}+\frac{1}{4}\right)^2}{4 t}} \sin (\pi  {x_1} {x_2}), \quad (x_1,x_2,t) \in Q = \Omega \times (0,T),
\end{equation}
defining the right-hand side $f$ and the inhomogeneous Dirichlet data on $\Sigma$. The inhomogeneous Dirichlet boundary condition is treated via homogenization as for the elliptic case, see \cite[p. 246]{Steinbach2008}. The spatial domain $\Omega$ is decomposed into uniform triangles with the uniform mesh size $h_x$ as given in Figure~\ref{Num:Fig:NetzeDreieckGlm} for level 0.

The temporal domain $(0,1/2) = (0,T)$ is decomposed into nonuniform elements with the nodes
\begin{equation} \label{Num:Zeinetz}
	t_0 = 0.0, \quad t_1 = 1/32, \quad t_2 = 1/16, \quad t_3 = 1/8, \quad t_4 = 1/2 = T.
\end{equation}

The assembling of the matrices  $A_{h_t}^{\mathcal H_T}$, $M_{h_t}^{\mathcal H_T}$, and $C_{h_t}^{\mathcal H_T}$ is done as proposed in \cite{Zank2020Exact}. The integrals for computing the projection $Q_h^0 f$ in \eqref{VF:Projektion} are calculated by using high-order quadrature rules. The solution $\vec u$ of the global linear system \eqref{VF:LGS} is solved in MATLAB by using the Bartels-Stewart methods (Algorithm~\ref{Loeser:BSr:Alg} and Algorithm~\ref{Loeser:BSc:Alg}) and the Fast Diagonalization method (Algorithm~\ref{Loeser:FD:Alg}), where the occurring spatial linear systems in Algorithm~\ref{Loeser:BSr:Alg}, Algorithm~\ref{Loeser:BSc:Alg} and Algorithm~\ref{Loeser:FD:Alg} are solved with the help of the sparse direct solver MUMPS 5.3.3 \cite{MUMPS2, MUMPS1} in the standard configuration. The other steps of Algorithm~\ref{Loeser:BSr:Alg}, Algorithm~\ref{Loeser:BSc:Alg} and Algorithm~\ref{Loeser:FD:Alg}, i.e. the real-Schur, complex-Schur, eigenvalue, singular value decompositions, and applying the transformation matrices, are realized by MATLAB routines. All calculations presented in this section were performed on a PC with two Intel Xeon CPUs E5-2687W v4 @ 3.00GHz, i.e. in sum 24 cores, and 512 GB main memory.

\subsection{Numerical Example for the Real Part of the Eigenvalues $\lambda_\ell$ and the Condition of the Transformation Matrix $X_t$} \label{Sec:Num:EW}

In this subsection, we investigate the eigenvalues $\lambda_\ell$ of $(A_{h_t}^{\mathcal H_T})^{-1} M_{h_t}^{\mathcal H_T}$ and the condition number of the transformation matrix $X_t$, occurring in the Fast Diagonalization method in Subsection~\ref{Sec:Loeser:FD}, of the corresponding eigenvectors of $(A_{h_t}^{\mathcal H_T})^{-1} M_{h_t}^{\mathcal H_T}.$ In Table~\ref{Num:EW:K2_X_t}, the smallest real part $\min_{\ell=1,\dots,N_t} \Re (\lambda_\ell)$ of the complex eigenvalues $\lambda_\ell$ of the eigenvalue decomposition \eqref{Loeser:EWP:EWP}, the minimal singular value $\sigma_{\min}(X_t)$, the maximal singular value $\sigma_{\max}(X_t)$, and the spectral condition number $\kappa_2(X_t)$ of the transformation matrix $X_t \in \C^{N_t \times N_t}$ of the eigenvalue decomposition \eqref{Loeser:FD:Zerlegung} are given for the nonuniform time mesh \eqref{Num:Zeinetz} with a uniform refinement strategy. The smallest real part $\min_{\ell=1,\dots,N_t} \Re (\lambda_\ell)$ is small but still strictly positive, see \eqref{Loeser:EWP:positiverRT}. The spectral condition number $\kappa_2(X_t)$ grows fast, which leads to numerical instability. However, the additional singular value decomposition \eqref{Loeser:FD:SVD} damps these instabilities, see Table~\ref{Num:FD:Tab:Fehler}. Further investigations of this issue are needed and will be done elsewhere.

\begin{table}[tbhp]
\begin{center}
\begin{tabular}{rcccccc} 
	\hline
 $N_t$ & $h_{t,\max}$ & $h_{t,\min}$ & $\min_{\ell=1,\dots,N_t} \Re (\lambda_\ell)$ & $\sigma_{\min}(X_t)$  & $\sigma_{\max}(X_t)$ & $\kappa_2(X_t)$ \\
		\hline   
   4 & 0.37500 & 0.03125 & 1.514e-02  & 2.041e-01  & 1.954e+00  & 9.576e+00  \\ 
   8 & 0.18750 & 0.01562 & 4.991e-03  & 4.049e-02  & 3.109e+00  & 7.678e+01  \\ 
  16 & 0.09375 & 0.00781 & 1.727e-03  & 2.174e-03  & 4.235e+00  & 1.948e+03  \\ 
  32 & 0.04688 & 0.00391 & 5.529e-04  & 1.566e-04  & 5.978e+00  & 3.816e+04  \\ 
  64 & 0.02344 & 0.00195 & 1.735e-04  & 1.377e-05  & 8.936e+00  & 6.488e+05  \\ 
 128 & 0.01172 & 0.00098 & 5.241e-05  & 1.416e-06  & 1.301e+01  & 9.187e+06  \\ 
 256 & 0.00586 & 0.00049 & 1.540e-05  & 1.640e-07  & 1.966e+01  & 1.199e+08  \\ 
 512 & 0.00293 & 0.00024 & 3.769e-06  & 1.705e-08  & 2.827e+01  & 1.658e+09  \\ 
1024 & 0.00146 & 0.00012 & 7.281e-07  & 4.131e-09  & 3.812e+01  & 9.229e+09  \\   
    \hline
  \end{tabular}
    \caption{Numerical results for the smallest real part $\min_{\ell=1,\dots,N_t} \Re (\lambda_\ell)$ of the eigenvalue decomposition \eqref{Loeser:EWP:EWP} and for the condition of the transformation matrix $X_t$ of the eigenvalue decomposition \eqref{Loeser:FD:Zerlegung} for $T=\frac{1}{2}$ for a uniform refinement strategy.} \label{Num:EW:K2_X_t}
\end{center}
\end{table}

\subsection{Bartels-Stewart Method with Real-Schur Decomposition} \label{Sec:Num:BSr}

This subsection deals with a numerical example for the Bartels-Stewart method with real-Schur decomposition, developed in Subsection~\ref{Sec:Loeser:BSr}, i.e. Algorithm~\ref{Loeser:BSr:Alg}. We consider the setting, which is described at the beginning of this section. In addition to this situation, the ordering and symbolic factorization steps of the sparse direct solver MUMPS 5.3.3 \cite{MUMPS2, MUMPS1} are performed only once, since the sparsity patterns of the system matrices in step 3 of Algorithm~\ref{Loeser:BSr:Alg} remain the same for $k=1,\dots,N_t$.

In Table~\ref{Num:BSr:Tab:Fehler}, the numerical results for the smooth solution $u$ in \eqref{Num:ExakteLsg}, when a uniform refinement strategy is applied as in Figure~\ref{Num:Fig:NetzeDreieckGlm}, are given, where unconditional stability is observed and the convergence rates in $\| \cdot \|_{L^2(Q)}$ and $| \cdot |_{H^1(Q)}$ are as expected from the error estimates \eqref{VF:Fehler_L2} and \eqref{VF:Fehler_H1}. The last column of Table~\ref{Num:BSr:Tab:Fehler} states the computation times in seconds of the Bartels-Stewart method with real-Schur decomposition (Algorithm~\ref{Loeser:BSr:Alg}), where the computing time for assembling the matrices $A_{h_t}^{\mathcal H_T}, M_{h_t}^{\mathcal H_T}, A_{h_x}, M_{h_x}$, and $C_{h_t}^{\mathcal H_T}$ is not included. We observe that the calculating time in Table~\ref{Num:BSr:Tab:Fehler} grows with factors 11.3, 9.6, 10.9 for the last three levels, which are smaller than the factor 16 resulting from the complexity $\mathcal O(\mathrm{dof}^{4/3})$ in Table~\ref{Loeser:BSr:Tab:Komplexitaet}.

\begin{table}[tbhp]
\begin{center}
\begin{tabular}{rccccccc|r} 
	\hline
  dof & $h_x$ & $h_{t,\max}$ & $h_{t,\min}$ & $\norm{u - \widetilde u_{h}}_{L^2(Q)}$  & eoc & $\abs{u - \widetilde u_{h}}_{H^1(Q)}$ & eoc & Solving\\
		\hline    
       20 & 0.354 & 0.375 & 0.0313 & 3.326e-01 & 0.00 & 4.314e+00 & 0.0  &    $\approx$ 0.0 \\ 
      264 & 0.177 & 0.188 & 0.0156 & 1.089e-01 & 1.30 & 2.702e+00 & 0.5  &    $\approx$ 0.0 \\ 
     2576 & 0.088 & 0.094 & 0.0078 & 3.136e-02 & 1.64 & 1.440e+00 & 0.8  &    $\approx$ 0.0 \\ 
    22560 & 0.044 & 0.047 & 0.0039 & 8.309e-03 & 1.84 & 6.984e-01 & 1.0  &     0.1 \\  
   188480 & 0.022 & 0.023 & 0.0020 & 2.127e-03 & 1.93 & 3.447e-01 & 1.0  &     0.6 \\ 
  1540224 & 0.011 & 0.012 & 0.0010 & 5.376e-04 & 1.96 & 1.707e-01 & 1.0  &     5.7 \\ 
 12452096 & 0.006 & 0.006 & 0.0005 & 1.352e-04 & 1.98 & 8.502e-02 & 1.0  &    64.3 \\ 
100139520 & 0.003 & 0.003 & 0.0002 & 3.393e-05 & 1.99 & 4.244e-02 & 1.0  &   615.7 \\ 
803210240 & 0.001 & 0.001 & 0.0001 & 8.500e-06 & 1.99 & 2.120e-02 & 1.0  &  6681.0 \\ 
    \hline
  \end{tabular}
    \caption{Numerical results of the Galerkin finite element discretization \eqref{VF:VF_HT_disk_Tensor_gestoert} for the L-shape \eqref{Num:Lshape} and $T=\frac{1}{2}$ for the function $u$ in \eqref{Num:ExakteLsg} for a uniform refinement strategy solved by the Bartels-Stewart method with real-Schur decomposition (Algorithm~\ref{Loeser:BSr:Alg}) using the sparse direct solver MUMPS 5.3.3, where the computation times are given in seconds.} \label{Num:BSr:Tab:Fehler}
\end{center}
\end{table}

\subsection{Bartels-Stewart Method with Complex-Schur Decomposition} \label{Sec:Num:BSc}

In this subsection, a numerical example for the Bartels-Stewart method with complex-Schur decomposition, developed in Subsection~\ref{Sec:Loeser:BSc}, i.e. Algorithm~\ref{Loeser:BSc:Alg}, is investigated. We consider the setting, which is described at the beginning of this section. In addition to this situation, the ordering and symbolic factorization steps of the sparse direct solver MUMPS 5.3.3 \cite{MUMPS2, MUMPS1} are performed only once, since the sparsity patterns of the system matrices in step 3 of Algorithm~\ref{Loeser:BSc:Alg} remain the same for $k=1,\dots,N_t$.

In Table~\ref{Num:BSc:Tab:Fehler}, the numerical results for the smooth solution $u$ in \eqref{Num:ExakteLsg}, when a uniform refinement strategy is applied as in Figure~\ref{Num:Fig:NetzeDreieckGlm}, are given, where errors and convergence rates are the same as for the Bartels-Stewart method with real-Schur decomposition (Table~\ref{Num:BSr:Tab:Fehler}). The last column of Table~\ref{Num:BSc:Tab:Fehler} states the computation times in seconds of the Bartels-Stewart method with complex-Schur decomposition (Algorithm~\ref{Loeser:BSc:Alg}), 
where the computing time for assembling the matrices $A_{h_t}^{\mathcal H_T}, M_{h_t}^{\mathcal H_T}, A_{h_x}, M_{h_x}$, and $C_{h_t}^{\mathcal H_T}$ is not included. We observe that the calculating time in Table~\ref{Num:BSc:Tab:Fehler} grows with factors 9.3, 10.7, 12.2 for the last three levels, which are smaller than the factor 16 resulting from the complexity $\mathcal O(\mathrm{dof}^{4/3})$ in Table~\ref{Loeser:BSr:Tab:Komplexitaet}. Moreover, we see that the Bartels-Stewart method with complex-Schur decomposition is slower than the real version, see Table~\ref{Num:BSr:Tab:Fehler}.

\begin{table}[tbhp]
\begin{center}
\begin{tabular}{rccccccc|r} 
	\hline
  dof & $h_x$ & $h_{t,\max}$ & $h_{t,\min}$ & $\norm{u - \widetilde u_{h}}_{L^2(Q)}$  & eoc & $\abs{u - \widetilde u_{h}}_{H^1(Q)}$ & eoc & Solving\\
		\hline    
       20 & 0.354 & 0.375 & 0.0313 & 3.326e-01 & 0.00 & 4.314e+00 & 0.0  &    $\approx$ 0.0 \\ 
      264 & 0.177 & 0.188 & 0.0156 & 1.089e-01 & 1.30 & 2.702e+00 & 0.5  &    $\approx$ 0.0 \\ 
     2576 & 0.088 & 0.094 & 0.0078 & 3.136e-02 & 1.64 & 1.440e+00 & 0.8  &    $\approx$ 0.0 \\ 
    22560 & 0.044 & 0.047 & 0.0039 & 8.309e-03 & 1.84 & 6.984e-01 & 1.0  &     0.4 \\  
   188480 & 0.022 & 0.023 & 0.0020 & 2.127e-03 & 1.93 & 3.447e-01 & 1.0  &     0.9 \\ 
  1540224 & 0.011 & 0.012 & 0.0010 & 5.376e-04 & 1.96 & 1.707e-01 & 1.0  &     9.7 \\ 
 12452096 & 0.006 & 0.006 & 0.0005 & 1.352e-04 & 1.98 & 8.502e-02 & 1.0  &    90.6 \\ 
100139520 & 0.003 & 0.003 & 0.0002 & 3.393e-05 & 1.99 & 4.244e-02 & 1.0  &   975.0 \\ 
803210240 & 0.001 & 0.001 & 0.0001 & 8.500e-06 & 1.99 & 2.120e-02 & 1.0  &  11872.6 \\ 
    \hline
  \end{tabular}
    \caption{Numerical results of the Galerkin finite element discretization \eqref{VF:VF_HT_disk_Tensor_gestoert} for the L-shape \eqref{Num:Lshape} and $T=\frac{1}{2}$ for the function $u$ in \eqref{Num:ExakteLsg} for a uniform refinement strategy solved by the Bartels-Stewart method with complex-Schur decomposition (Algorithm~\ref{Loeser:BSc:Alg}) using the sparse direct solver MUMPS 5.3.3, where the computation times are given in seconds.} \label{Num:BSc:Tab:Fehler}
\end{center}
\end{table}

\subsection{Fast Diagonalization Method} \label{Sec:Num:FD}

In this subsection, a numerical example for the Fast Diagonalization method, developed in Subsection~\ref{Sec:Loeser:FD}, i.e. Algorithm~\ref{Loeser:FD:Alg}, is given. We consider the setting, which is described at the beginning of this section. In addition to this situation, time parallelization, but no spatial parallelization is applied, i.e. the $N_t$ spatial problems of step 3 in Algorithm~\ref{Loeser:FD:Alg} can be solved in parallel if $N_t$ cores are available.

In Table~\ref{Num:FD:Tab:Fehler}, the numerical results for the smooth solution $u$ in \eqref{Num:ExakteLsg}, when a uniform refinement strategy is applied as in Figure~\ref{Num:Fig:NetzeDreieckGlm}, are given, where errors and convergence rates are the same as for the Bartels-Stewart methods (Tables~\ref{Num:BSr:Tab:Fehler} and~\ref{Num:BSc:Tab:Fehler}). For the last level in Table~\ref{Num:FD:Tab:Fehler}, the $L^2$ error is slightly larger than the corresponding error for the Bartels-Stewart methods in Table~\ref{Num:BSr:Tab:Fehler} and Table~\ref{Num:BSc:Tab:Fehler}, due to the large condition number of the transformation matrix $X_t$, see Table~\ref{Num:EW:K2_X_t}. The last column of Table~\ref{Num:FD:Tab:Fehler} states the computation times in seconds of the Fast Diagonalization method (Algorithm~\ref{Loeser:FD:Alg}), where the computing time for assembling the matrices $A_{h_t}^{\mathcal H_T}, M_{h_t}^{\mathcal H_T}, A_{h_x}, M_{h_x}$, and $C_{h_t}^{\mathcal H_T}$ is not included. We observe that the calculating time in Table~\ref{Num:FD:Tab:Fehler} grows with factors 8.2, 8.8, 10.1 for the last three levels, which are smaller than the factor 16 resulting from the complexity $\mathcal O(\mathrm{dof}^{4/3})$ in Table~\ref{Loeser:BSr:Tab:Komplexitaet}. Additionally, we see that the Fast Diagonalization method is much faster than the Bartels-Stewart methods (Table~\ref{Num:BSr:Tab:Fehler}, Table~\ref{Num:BSc:Tab:Fehler}) due to the time parallelization.

\begin{table}[tbhp]
\begin{center}
\begin{tabular}{rccccccc|r} 
	\hline
  dof & $h_x$ & $h_{t,\max}$ & $h_{t,\min}$ & $\norm{u - \widetilde u_{h}}_{L^2(Q)}$  & eoc & $\abs{u - \widetilde u_{h}}_{H^1(Q)}$ & eoc & Solving\\
		\hline    
       20 & 0.354 & 0.375 & 0.0313 & 3.326e-01 & 0.00 & 4.314e+00 & 0.0  &    $\approx$ 0.0 \\ 
      264 & 0.177 & 0.188 & 0.0156 & 1.089e-01 & 1.30 & 2.702e+00 & 0.5  &    $\approx$ 0.0 \\ 
     2576 & 0.088 & 0.094 & 0.0078 & 3.136e-02 & 1.64 & 1.440e+00 & 0.8  &    $\approx$ 0.0 \\ 
    22560 & 0.044 & 0.047 & 0.0039 & 8.309e-03 & 1.84 & 6.984e-01 & 1.0  &     0.1 \\  
   188480 & 0.022 & 0.023 & 0.0020 & 2.127e-03 & 1.93 & 3.447e-01 & 1.0  &     0.3 \\ 
  1540224 & 0.011 & 0.012 & 0.0010 & 5.376e-04 & 1.96 & 1.707e-01 & 1.0  &     0.9 \\ 
 12452096 & 0.006 & 0.006 & 0.0005 & 1.352e-04 & 1.98 & 8.502e-02 & 1.0  &    7.4 \\ 
100139520 & 0.003 & 0.003 & 0.0002 & 3.393e-05 & 1.99 & 4.244e-02 & 1.0  &   64.9 \\ 
803210240 & 0.001 & 0.001 & 0.0001 & 8.855e-06 & 1.94 & 2.121e-02 & 1.0  &  652.8 \\ 
    \hline
  \end{tabular}
    \caption{Numerical results of the Galerkin finite element discretization \eqref{VF:VF_HT_disk_Tensor_gestoert} for the L-shape \eqref{Num:Lshape} and $T=\frac{1}{2}$ for the function $u$ in \eqref{Num:ExakteLsg} for a uniform refinement strategy solved by the Fast Diagonalization method (Algorithm~\ref{Loeser:FD:Alg}) using the sparse direct solver MUMPS 5.3.3, where the computation times are given in seconds.} \label{Num:FD:Tab:Fehler}
\end{center}
\end{table}

\section{Conclusions} \label{Sec:Zum}

In this work, we studied efficient direct solvers for the global linear system arising from the space-time Galerkin finite element discretization of parabolic initial-boundary value problems in anisotropic Sobolev spaces in combination with the Hilbert-type transformation operator $\mathcal H_T$. Two algorithms based on the Bartels-Stewart method and one algorithm based on the Fast Diagonalization method were developed and analyzed. The latter allows a complete parallelization in time. We gave complexity estimates for these three algorithms. We presented numerical experiments for a two-dimensional spatial domain, where the spatial subproblems were solved by sparse direct solvers. These numerical results confirmed the efficient applicability of the space-time approach in anisotropic Sobolev spaces in connection with the direct space-time solvers proposed in the paper.

For the spatial subproblems occurring in the algorithms, (preconditioned) iterative solvers  can also  be used. We only used piecewise linear ansatz and test functions, but the approach can easily be generalized to shape functions of an arbitrary polynomial degree, and to graded or even adaptive meshes in space. Furthermore, the space-time approach presented in this paper also works for autonomous parabolic problems with diffusion coefficients depending on the spatial variable only, and even for non-autonomous parabolic problems with diffusion coefficients being a product of a function in $x$ and a function in $t$. Moreover, the direct solvers proposed in this paper can be used in connection with a preconditioner for more general linear and even non-linear parabolic problems.

\bibliographystyle{acm}
\bibliography{literatur}

\end{document}